\documentclass[12pt]{article}
\usepackage{hyperref}
\usepackage{amssymb,amsmath}
\usepackage{amsthm}
\newtheorem*{theorem}{Theorem}
\usepackage{color}
\usepackage{graphicx}
\usepackage[symbol]{footmisc}
\usepackage[ngerman,english]{babel}
\usepackage[dvipsnames]{xcolor}
\usepackage[top=2.5cm,bottom=2.5cm,left=2.5cm,right=2.5cm]{geometry}

\begin{document}
\begin{center}
\LARGE{\textbf{A root finding method with
\\arbitrary order of convergence }}
$$
$$
\large
{Alois Schiessl}
\footnote[1]{E-Mail: \texttt{aloisschiessl@posteo.de}}
\vskip 0.5cm
\centerline{\today}
\end{center}
\[\]
\[\]
\[\]
\begin{abstract}
Let $a\in \mathbb{R}^{+}\backslash\left\{0\right\}$ and $M\in\mathbb{N}$. We consider the equation $t^M-a=0$, which is equivalent to $1-\frac{t^M}{a}=0\,.$ The real solution is $\sqrt[M]{a}$. In this publication, we present a method that enables the calculation of $\sqrt[M]{a}$ with arbitrary order of convergence using only polynomials. We define the fixed point function
\[
F\left(x\right)
=\prod_{\ell=1}^{P}\left(1+\frac{1}{\ell\cdot M}\right)
\int\limits_{0}^{x}\!\left(1-{\frac{{t}^{M}}{a}}\right)^{P}{\rm d}t
=\sum\limits_{k=0}^{P}\frac{\left(-1\right)^{\,k}}{a^{\,k}}\cdot\binom{P}{k}\cdot\frac{x^{\,k\,\cdot M+1}}{k\,\cdot M+1}
\]
This is a polynomial of degree $\left(P\cdot M+1\right)$ with $\left(P+1\right)$ terms. The calculation of $\sqrt[M]{a}$ is thus reduced to a polynomial evaluation. The computational tests we performed demonstrate the efficiency of the method.
\[\]
\[\textbf{Zusammenfassung}\]
Es sei $a\in \mathbb{R}^{+}\backslash\left\{0\right\}$ und $M\in\mathbb{N}$. Vorgelegt ist die Gleichung $t^M-a=0$, die äquivalent zu $1-\frac{t^M}{a}=0$ ist. Die reelle L\"osung hiervon ist $\sqrt[M]{a}$. In dieser Ver\"offentlichung stellen wir ein Verfahren vor, das die Berechnung von $\sqrt[M]{a}$ mit beliebiger Konvergenzordnung erm\"oglicht und nur Polynome verwendet. Wir definieren die Fixpunktfunktion
\[F\left(x\right)
=\prod_{\ell=1}^{P}\left(1+\frac{1}{\ell\cdot M}\right)
\int\limits_{0}^{x}\!\left(1-{\frac{{t}^{M}}{a}}\right)^{P}{\rm d}t
=\sum\limits_{k=0}^{P}\frac{\left(-1\right)^{\,k}}{a^{\,k}}\cdot\binom{P}{k}\cdot\frac{x^{\,k\,\cdot M+1}}{k\,\cdot M+1}
\]
Das ist ein Polynom vom Grad $\left(P\cdot M+1\right)$ mit $\left(P+1\right)$ Summanden. Die Berechnung von $\sqrt[M]{a}$ reduziert sich dadurch zu einer Polynomauswertung. Anhand ausgew\"ahlter Beispiele von Wurzelberechnungen zeigen wir die Effizienz des Verfahrens.
\[\text{Deutsche Version ab Seite 19}\]
\[\]
\end{abstract}

\selectlanguage{english}

\section{Introduction and Main Result}
One of the best-known method for calculating an approximate value for $\sqrt[M]a$ for a positive real number $a$ and $M\in \mathbb{N}$ is to apply Newton's method \cite{WG}
\[
x_{n+1}=x_n-\frac{f\left(x_n\right)}{f\,'\left(x_n\right)}
\]
to the function
\[f\left(x\right)=x^{M}-a\;,\]
from which the fixed point iteration can be derived
\[x_{n+1}
=x_n-\frac{x_n^{M}-a}{M\,x_n^{M -1}}
=x_n-\frac{1}{M}\left(x_n-\frac{a}{x_n^{M -1}}\right)\;;\;n=0,1,\ldots\]
If the initial value $x_0$ is sufficiently close to $\sqrt[M]a$, the sequence converges quadratically to $\sqrt[M]a$. Each additional Iteration approximately doubles the number of correct digits. In this article, we present a method that only works with polynomials and allows for any order of convergence. The method is based on the function
\[f\left(t\right)=1-\frac{t^M}{a}\;.\]
Now we state our main result:
\begin{theorem}
$ \\ $
Let $a\in \mathbb{R}\backslash\left\{0\right\}$ and $M\in \mathbb{N}$. We define
\[f: \mathbb{R} \to \mathbb{R}\;;\]
\[t\mapsto f\left(t\right)\;;\]
\[f\left(t\right)=1-\frac{t^M}{a}\;.\]
Let $P\in\mathbb{N}.$
We define
\[F: \mathbb{R} \to \mathbb{R}\;;\]
\[x\mapsto F\left(x\right)\;;\]
\[F\left(x\right)
=\prod_{\ell=1}^{P}\left(1+\frac{1}{\ell\cdot M}\right)\int\limits_{0}^{x}\!f\left(t\right)^{P}{\rm d}t
=\prod_{\ell=1}^{P}\left(1+\frac{1}{\ell\cdot M}\right)
\int\limits_{0}^{x}\!\left(1-{\frac{{t}^{M}}{a}}\right)^{P}{\rm d}t\;.
\]
Let $a>0.$ Then, there hold the following statements for $F\left(x\right)$ 
\begin{itemize}
\item[(a)]
$\sqrt[M]{a}$ is a fixed point, that means: $F\left(\sqrt[M]{a}\right)=\sqrt[M]{a}\,.$
\item[(b)]
The derivatives up to $\left(P+1\right)$ at the fixed point $\sqrt[M]{a}$ are given by
\begin{align*}
F{\left(\sqrt[M]{a}\right)}^{\left(k\right)}=\left\{
\begin{aligned}
0\qquad\qquad\qquad\qquad 1\,\le &\,k\,\leq{P}
\\
\\
\frac{\left(-1\right)^P}{\sqrt[M]{a^P}}\prod_{\ell=1}^{P}\left(1+\ell\cdot M\right)
\qquad\qquad\quad k=&P+1
\\ 
\end{aligned} 
\right.
\end{align*}
\item[(c)]
Let $x_0$ an initial value sufficiently close to $\sqrt[M]{a}\,,$ then the sequence
\[
x_{n+1}=F\left(x_{{n}}\right)\;;\;n=0,1,\ldots
\]
converges to $\sqrt[M]{a}$ with order of convergence exactly $\left(P+1\right).$
\end{itemize}
\end{theorem}

\section{Proof of the theorem}
In this section we give a proof of the main result using only Classical Analysis.
\subsection{Fixed point}
We have to prove that for $a>0$ the following statement is true
\[F\left(\sqrt[M]{a}\right)=\prod_{\ell=1}^{P}\left(1+\frac{1}{\ell\cdot M}\right)
\int\limits_{0}^{\sqrt[M]{a}}\!\left(1-{\frac{{t}^{M}}{a}}\right)^{P}{\rm d}t=\sqrt[M]{a}\;.\]
We will first consider the definite integral
\[\int\limits_{0}^{\sqrt[M]{a}}\!\left(1-\frac{{t}^{M}}{a}\right)^{P}{\rm d}t\;.\]
We make the substitution $\frac{t^M}{a}=u$ to obtain 
\[\frac{{\rm d} \,t} {{\rm d}\, u}=\sqrt[M]{a}\cdot\frac{1}{M}\cdot u^{\frac{1}{M}-1}\;.\]
Therefore: 
\[
\int\limits_{0}^{\sqrt[M]{a}}\!\left(1-{\frac{{t}^{M}}{a}}\right)^{P}{\rm d}t
=\int\limits_{0}^{1}\!\left(1-u\right)^{P}\sqrt[M]{a}\cdot\frac{1}{M}\cdot u^{\frac{1}{M}-1}{\rm d}u
=\frac{\sqrt[M]{a}}{M}\int\limits_{0}^{1}\!u^{\frac{1}{M}-1}\cdot\left(1-u\right)^{P}{\rm d}u\;.
\]
Since the lower limit $t=0$ was replaced with $u=0$, and the upper limit $\sqrt[M]{a}$ with $1$. The next task is now to evaluate the integral
\[
\int\limits_{0}^{1}\!u^{\frac{1}{M}-1}\cdot\left(1-u\right)^{P}{\rm d}u\;.
\]
Let us remind the real Beta function \cite{AS}:
\[{\rm B}\left(\alpha,\beta\right)=\int\limits_{0}^{1}\!{u}^{\alpha-1}\left(1-u\right)^{\beta-1}\,{\rm d}u
\;;\,\alpha>0\;;\,\beta>0\;,\]
which is closely related to the gamma function \cite{EA}
\[{\rm B}\left(\alpha,\beta\right)
={\frac{\Gamma\left(\alpha\right)\Gamma\left(\beta\right)}{\Gamma\left(\alpha+\beta\right)}}\;.\]
Substituting
\[\alpha=\frac{1}{M}\;;\quad\beta=P+1\]
gives a representation by the $\Gamma$ function
\[
\int\limits_{0}^{1}\!u^{\frac{1}{M}-1}\cdot\left(1-u\right)^{P}{\rm d}u
=\frac{\Gamma\left(\frac{1}{M}\right)\cdot\Gamma\left(P+1\right)}{\Gamma\left(P+1+\frac{1}{M}\right)}\;.
\]
Now, $P$ is a natural number, so $\Gamma\left(P+1\right)=P\,!\;$. This simplifies the term to
\[
\frac{\Gamma\left(\frac{1}{M}\right)\cdot\Gamma\left(P+1\right)}{\Gamma\left(P+1+\frac{1}{M}\right)}
=\frac{\Gamma\left(\frac{1}{M}\right)P\,!}{\Gamma\left(P+1+\frac{1}{M}\right)}\;.
\]
Next, we focus on
\[
\Gamma\left(P+1+\frac{1}{M}\right)\;.
\]
Repeated using of the definition of the $\Gamma$-Function
\[\Gamma\!\left(1+x\right)=x\,\Gamma\!\left(x\right)\]
with $x=\frac{1}{M}$ yields step by step
\[\Gamma\!\left(1+\frac{1}{M}\right)=\frac{\Gamma\!\left(\frac{1}{M}\right)}{M}\;;\]
\[
\Gamma\!\left(2+\frac{1}{M}\right)
=\left(1+\frac{1}{M}\right)\frac{\Gamma\!\left(\frac{1}{M}\right)}{M}
=\prod\limits_{\ell=1}^{1}\left(1+\frac{1}{\ell\cdot M}\right)\frac{1\,!\cdot\Gamma\left(\frac{1}{M}\right)}{M}\;;
\]
\[
\Gamma\!\left(3+\frac{1}{M}\right)
=\left(1+\frac{1}{M}\right)\left(2+\frac{1}{M}\right)\frac{\Gamma\!\left(\frac{1}{M}\right)}{M}
=\prod\limits_{\ell=1}^{2}\left(1+\frac{1}{\ell\cdot M}\right)\frac{2\,!\cdot\Gamma\left(\frac{1}{M}\right)}{M} \;;
\]
\[
\Gamma\!\left(4+\frac{1}{M}\right)
=\left(1+\frac{1}{M}\right)\left(2+\frac{1}{M}\right)\left(3+\frac{1}{M}\right)\frac{\Gamma\!\left(\frac{1}{M}\right)}{M}\\
=
\prod\limits_{\ell=1}^{3}\left(1+\frac{1}{\ell\cdot M}\right)\frac{3\,!\cdot\Gamma\left(\frac{1}{M}\right)}{M}\;;
\]
\[\vdots\]
By using an induction argument it is a simple matter to prove that
\[
\Gamma\left(P+1+\frac{1}{M}\right)
=\prod\limits_{\ell=1}^{P}\left(1+\frac{1}{\ell\cdot M}\right)\frac{P\,!\cdot\Gamma\left(\frac{1}{M}\right)}{M}\;.
\]
Thus, for the definite integral, we obtain
\[
\int\limits_{0}^{1}\!u^{\frac{1}{M}-1}\cdot\left(1-u\right)^{P}{\rm d}u
=\frac{\Gamma\left(\frac{1}{M}\right)\cdot\Gamma\left(P+1\right)}{\Gamma\left(P+1+\frac{1}{M}\right)}
=\frac{\Gamma\left(\frac{1}{M}\right)\cdot P\,!}{\prod\limits_{\ell=1}^{P}\left(1+\frac{1}{\ell\cdot M}\right)\frac{P\,!\cdot\Gamma\left(\frac{1}{M}\right)}{M}}
=\frac{M}{\prod\limits_{\ell=1}^{P}\left(1+\frac{1}{\ell\cdot M}\right)}\;.
\]
Now we just need to put everything together. That gives us
\begin{align*}
F\left(\sqrt[M]{a}\right)
&=\prod_{\ell=1}^{P}\left(1+\frac{1}{\ell\cdot M}\right)\int\limits_{0}^{\sqrt[M]{a}}\!\left(1-{\frac{{t}^{M}}{a}}\right)^{P}{\rm d}t
\\
&=\prod_{\ell=1}^{P}\left(1+\frac{1}{\ell\cdot M}\right)\frac{\sqrt[M]{a}}{M}\int\limits_{0}^{1}\!u^{\frac{1}{M}-1}\cdot\left(1-u\right)^{P}{\rm d}u
\\
&=\prod_{\ell=1}^{P}\left(1+\frac{1}{\ell\cdot M}\right)\frac{\sqrt[M]{a}}{M}\cdot \frac{M}{\prod\limits_{\ell=1}^{P}\left(1+\frac{1}{\ell\cdot M}\right)}
\\
&=\prod_{\ell=1}^{P}\left(1+\frac{1}{\ell\cdot M}\right)\frac{\sqrt[M]{a}}{\prod\limits_{\ell=1}^{P}\left(1+\frac{1}{\ell\cdot M}\right)}
=\sqrt[M]{a}\;.
\end{align*}
The proof of the fixed point statement is completed.
\subsection{Derivatives}
Next, we focus on the derivatives.
\subsubsection{First Derivative}
Starting with the first derivative we obtain immediately
\[
F\,'\left(x\right)
=\left(\prod_{\ell=1}^{P}\left(1+\frac{1}{\ell\cdot M}\right)\int\limits_{0}^{x}\!\left(1-{\frac{{t}^{M}}{a}}\right)^{P}{\rm d}t\right)'
=\prod_{\ell=1}^{P}\left(1+\frac{1}{\ell\cdot M}\right)\left(1-{\frac{{x}^{M}}{a}}\right)^{P}\;.
\]
Plugging in $\sqrt[M]{a}$ yields $F\,'\left(\sqrt[M]{a}\right)=0$. We will need this result later.
\subsubsection{Higher derivatives}
So far, we have a formula for the first derivative
\[
F\,'\left(x\right)
=\prod_{\ell=1}^{P}\left(1+\frac{1}{\ell\cdot M}\right)\cdot\left(1-\frac{x^M}{a}\right)^{P}.
\]
We still have to prove that the following statement holds at the fixed point $x=\sqrt[M]{a}\,:$
\begin{align*}
F{\left(\sqrt[M]{a}\right)}^{\left(k\right)}=\left\{
\begin{aligned}
0\qquad\qquad\qquad\qquad 1\,\le &\,k\,\leq{P}
\\
\\
\frac{\left(-1\right)^P}{\sqrt[M]{a^P}}\prod_{\ell=1}^{P}\left(1+\ell\cdot M\right)
\qquad\qquad\quad k=&P+1
\\ 
\end{aligned} 
\right.
\end{align*}
This requires a little bit more work. We must differentiate $\left(P+1\right)$ times $F\left(x\right)$ and then plugging in $x=\sqrt[M]{a}$. Using the first derivative $F\,'\left(x\right)$ the task simplifies to 
\begin{align*}
\left.F\left(x\right)^{\left(P+1\right)}\right|_{x={\sqrt[M]{a}}}=\left.\bigl(F\,'\left(x\right)\bigr)^{\left(P\right)}\right|_{x={\sqrt[M]{a}}}
=\left.\left(\prod_{\ell=1}^{P}\left(1+\frac{1}{\ell\cdot M}\right)\cdot\left(1-\frac{x^M}{a}\right)^{P}\right)^{\left(P\,\right)}\right|_{x={\sqrt[M]{a}}}
\end{align*}
On the right-hand side, we need at $x=\sqrt[M]{a}$ the $P^{\,th}$ derivative of
\[\left(1-\frac{x^M}{a}\right)^P\]
Let us put
\[f\left(x\right)=\left(1-\frac{x^M}{a}\right)\]
and first calculating some particular values of $f\left(x\right),$ needed later:
\[f\left(\sqrt[M]{a}\right)=0\;;\]
\[f'\left(\sqrt[M]{a}\right)=\left.\left(1-\frac{x^M}{a}\right)'\right|_{x={\sqrt[M]{a}}}
=\ldots=\frac{-M}{\sqrt[M]{a}}\;.\]
In the next step, we expand $f\left(x\right)$ into a Taylor series of degree 1 at the fixed point $x_0=\sqrt[M]{a}$ with Lagrange remainder:
\[
f\left(x\right)
=f\left(\sqrt[M]{a}\right)
+f'\left(\sqrt[M]{a}\right)\left(x-\sqrt[M]{a}\right)
+\frac{f''\left(\xi\right)}{2}\left(x-\sqrt[M]{a}\right)^{2}
\]
where $\xi$ is between $x$ and $\sqrt[M]{a}$.
$ \\ $
$ \\ $
Since $f\left(\sqrt[M]{a}\right)=0$ and $f'\left(\sqrt[M]{a}\right)=\frac{-M}{\sqrt[M]{a}}$, we get
\[
f\left(x\right)
=
\frac{-M}{\sqrt[M]{a}}\left(x-\sqrt[M]{a}\right)
+\frac{f''\left(\xi\right)}{2}\left(x-\sqrt[M]{a}\right)^{2}
\]
and can thus factor out $\left(x-\sqrt[M]{a}\right):$
\[
f\left(x\right)=
\left(x-\sqrt[M]{a}\right)
\left(\frac{-M}{\sqrt[M]{a}}+\frac{f''\left(\xi\right)}{2}\left(x-\sqrt[M]{a}\right)\right)\;.
\]
Next we raise $f\left(x\right)$ to the $P^{\,th}$ power\,:
\[
f\left(x\right)^{P}=
\left(x-\sqrt[M]{a}\right)^P
\left(\frac{-M}{\sqrt[M]{a}}+\frac{f''\left(\xi\right)}{2}\left(x-\sqrt[M]{a}\right)\right)^P\;.
\]
We calculate the higher derivatives of $f\left(x\right)^{P}$ by using the Leibniz rule \cite{HJ}.
Leibniz rule states that if $u\left(x\right)$ and $v\left(x\right)$ are $n$-times differentiable functions, then the product $u\left(x\right)\cdot v\left(x\right)$ is also $n$-times differentiable and its $n^{\,th}$ derivative is given by the sum
\[
\bigl(u\left(x\right)\cdot v\left(x\right)\bigr)^{\left(n\right)}
=\sum\limits_{k=0}^{n}\binom{n}{k}\,u\left(x\right)^{\left(k\right)}\cdot v\left(x\right)^{\left(n-k\right)}\;.
\]
Here we want to understand $u\left(x\right)^{\left(0\right)}=u\left(x\right)$ and $v\left(x\right)^{\left(0\right)}=v\left(x\right)\;.$
$\\ $
$\\ $
In our case we have:
\[
n=P\;;
\]
\[
u\left(x\right)=\left(x-\sqrt[M]{a}\right)^{P}\;;
\]
\[
v\left(x\right)=\left(\frac{-M}{\sqrt[M]{a}}+\frac{f''\left(x\right)}{2}\left(x-\sqrt[M]{a}\right)\right)^{P}\;.
\]
Thus we obtain
\[
f\left(x\right)^{P}=
\underbrace{{{\left(x-\sqrt[M]{a}\right)}^{P}}}_{u\left(x\right)}
\times
\underbrace{{{\left(\frac{-M}{\sqrt[M]{a}}+\frac{f''{{\left(\xi\right)}}}{2}\left(x-\sqrt[M]{a}\right)\right)}^{P}}}_{v\left(x\right)}\;.
\]
Now applying Leibniz rule
\[
\left(f\left(x\right)^{P}\right)^{\left({P}\right)}
=\sum\limits_{k=0}^{{P}}\binom{{P}}{k}
\,\left(\left(x-\sqrt[M]{a}\right)^{P}\right)^{\left(k\right)}
\times
\left(\left(\frac{-M}{\sqrt[M]{a}}+\frac{f''{\left(\xi\right)}}{2}\left(x-\sqrt[M]{a}\right)\right)^{{P}}\right)^{\left({P-k}\right)}\,.
\]
$ \\ $
We first focus on the sequence of derivatives
\[
u\left(x\right)^{\left(k\right)}
=\Bigl(\left(x-\sqrt[M]{a}\right)^{P}\Bigr)^{\left(k\right)},\;k=0,1,\ldots,P\;.
\]
From a collection of mathematical formulas \cite{formula}, we derive the following statement
\begin{align*}
\Bigl(\left(x-\sqrt[M]{a}\right)^{P}\Bigr)^{\left(k\right)}
=\frac{P\,!}{\left(P-k\right)\,!}\left(x-\sqrt[M]{a}\right)^{P-k}\;.
\end{align*}
$\\ $
In particular, we note that for $0\leq k <P$ in the formula
\[
u\left(x\right)^{\left(k\right)}=
\left(\left(x-\sqrt[M]{a}\right)^{{P}}\right)^{\left(k\right)}
=\frac{P\,!}{\left(P-k\right)\,!}\left(x-\sqrt[M]{a}\right)^{P-k},
\]
the factor
$\left(x-\sqrt[M]{a}\right)$ always occurs and thus for $x=\sqrt[M]{a}$ the terms vanish. Only the case of $k=P$ leads to a non-zero result, namely the product
\[
{\left.u{{\left(x\right)}^{\left(P\right)}}\right|}_{x=\sqrt[M]{a}}=\frac{P\,!}{\left(P-P\right)\,!}=P\,!\;.
\]
From the sum
\[
\left(f\left(x\right)^{P}\right)^{\left({P}\right)}
=\sum\limits_{k=0}^{{P}}\binom{{P}}{k}
\,\left(\left(x-\sqrt[M]{a}\right)^{P}\right)^{\left(k\right)}
\times
\left(\left(\frac{-M}{\sqrt[M]{a}}+\frac{f''{\left(\xi\right)}}{2}\left(x-\sqrt[M]{a}\right)\right)^{{P}}\right)^{\left({P-k}\right)}
\]
only the term with $k=P$ remains at the fixed point $x=\sqrt[M]{a}$. We can evaluate
\[
\left.\left(f\left(x\right)^{P}\right)^{\left({P}\right)}\right|_{x=\sqrt[M]{a}}=
\underbrace{\left.\left(\left(x-\sqrt[M]{a}\right)^{P}\right)^{\left(P\right)}\right|_{x=\sqrt[M]{a}}}_{=P\,!}
\times
\underbrace{\left.\left(\left(\frac{-M}{\sqrt[M]{a}}
+\frac{f''{\left(\xi\right)}}{2}\left(x-\sqrt[M]{a}\right)\right)^{{P}}\right)\right|_{x=\sqrt[M]{a}}}_{\left(\frac{-M}{\sqrt[M]{a}}\right)^P}
\]
\[=P\,!\cdot\left(\frac{-M}{\sqrt[M]{a}}\right)^P=P\,!\cdot\frac{\left(-1\right)^P M^P}{\sqrt[M]{a^P}}\;.\]
Our task was the calculation of
\begin{align*}
F\left(\sqrt[M]{a}\right)^{\left(P+1\right)}
&=\left.F\left(x\right)^{\left(P+1\right)}\right|_{x={\sqrt[M]{a}}}
=\left.F'\left(x\right)^{\left(P\right)}\right|_{x={\sqrt[M]{a}}}
\\
&=\left.\left(\prod_{\ell=1}^{P}\left(1+\frac{1}{\ell\cdot M}\right)\cdot\left(1-\frac{x^M}{a}\right)^{P}\right)^{\left(P\,\right)}\right|_{x={\sqrt[M]{a}}}\;.
\end{align*}
Now plugging in, we get
\begin{align*}
F\left(\sqrt[M]{a}\right)^{\left(P+1\right)}
=\prod_{\ell=1}^{P}\left(1+\frac{1}{\ell\cdot M}\right)\cdot P\,!\cdot\frac{\left(-1\right)^P \, M^P}{\sqrt[M]{a^P}}\;.
\end{align*}
Using the identity
\[
\prod_{\ell=1}^{P}\left(1+\frac{1}{\ell\cdot M}\right)
=\prod_{\ell=1}^{P}\left({\frac{\ell\cdot M+1}{\ell\cdot M}}\right)
=\frac{1}{P\,!\,M^P}\prod_{\ell=1}^{P}\left(1+\ell\cdot M\right)\,,
\]
the result can be reduced to a simpler form
\begin{align*}
F\left(\sqrt[M]{a}\right)^{\left(P+1\right)}
=\frac{1}{P!\,M^P}\prod_{\ell=1}^{P}\left(1+\ell\cdot M\right)\cdot\,P\,!\cdot\frac{\left(-1\right)^P M^P}{\sqrt[M]{a^P}}
=\frac{\left(-1\right)^P}{\sqrt[M]{a^P}}\,\prod_{\ell=1}^{P}\left(1+\ell\cdot M\right)\;.
\end{align*}
Putting all the pieces together gives us the desired result
\begin{align*}
F{\left(\sqrt[M]{a}\right)}^{\left(k\right)}=\left\{
\begin{aligned}
0\qquad\qquad\qquad\qquad 1\,\le &\,k\,\leq{P}
\\
\\
\frac{\left(-1\right)^P}{\sqrt[M]{a^P}}\prod_{\ell=1}^{P}\left(1+\ell\cdot M\right)
\qquad\qquad\quad k=&P+1
\\ 
\end{aligned} 
\right.
\end{align*}
This completes the proof of the derivations.

\subsection{Proof of convergence}
We do the prove by using Banach's fixed point theorem \cite{JM}, in a special form for continuously differentiable functions in $\mathbb{R}$.
$ \\ $
$ \\ $
\underline{Banach fixed point theorem} 
\\
Let $U\subseteq R$ be a closed subset; furthermore, let $F:U\rightarrow R$ be a mapping with the following properties\[
F\left(U\right)\subseteq U \;.
\]
There exists $0<L<1$ such that
\[
\left|F\,'\left(x\right)\right| \leq L\,,\; \text{for\,\;all} \;x\in U \;.
\]
Then the following statements hold:
\begin{itemize}
\item[(a)]
There is exactly one fixed point ${{x}^{*}}\in U$ of $F\;.$
\item[(b)]
For every initial value $x_0\in \mathbb{U}$, the sequence
$x_{n+1}=F\left(x_{{n}}\right)$
converges to $x^*\,.$
\end{itemize}
Now back to our concern. We have already calculated $F\left(\sqrt[M]{a}\right)=\sqrt[M]{a}$ and $F\,'\left(\sqrt[M]{a}\right)=0$. Then, for reasons of continuity, there exists a neighbourhood at the fixed point $\sqrt[M]{a}$
\[U_{\delta}=\left\{x\in\mathbb{R}:\left|x-\sqrt[M]{a}\right|\leq\delta\right\}\]
in which the following applies: $F\left(U\right)\subseteq U$ and $\left|F\,'\left(x\right)\right|<L<1$ for all $x\in U_{\delta}$.
Thus, the conditions for the fixed point theorem are satisfied and convergence follows.

\subsection{Order of convergence}
Finally, we need to proof that the order of convergence is exactly $\left(P+1\right)$.
We expand $F\left(x\right)$ into a Taylor series around the fixed point $x_0={\sqrt[M]{a}}\;:$
\begin{align*}
F\left(x\right)&=
F\left({\sqrt[M]{a}}\right)
+F'\left({\sqrt[M]{a}}\right)\left(x-{\sqrt[M]{a}}\right)
+\ldots
\\
&+F\left({\sqrt[M]{a}}\right)^{\left(P\right)}
\frac{\left(x-{\sqrt[M]{a}}\right)^{P}}{P\,!}
+F\left(\xi\right)^{\left(P+1\right)}
\frac{\left(x-{\sqrt[M]{a}}\right)^{P+1}}{\left(P+1\right)!}\;,
\end{align*}
where $\xi$ is between $x$ and ${\sqrt[M]{a}}$. Let us remind the fixed point property $F\left({\sqrt[M]{a}}\right)={\sqrt[M]{a}}$ and the differentiation results
\begin{align*}
F{\left(\sqrt[M]{a}\right)}^{\left(k\right)}=\left\{
\begin{aligned}
0\qquad\qquad\qquad\qquad 1\,\le &\,k\,\leq{P}
\\
\\
\frac{\left(-1\right)^P}{\sqrt[M]{a^P}}\prod_{k=1}^{P}\left(k\,\cdot M\right)
\qquad\qquad\quad k=&P+1
\\ 
\end{aligned} 
\right.
\end{align*}
The Taylor series
\begin{align*}
F\left(x\right)
&=\underbrace{F\left({\sqrt[M]{a}}\right)}_{={\sqrt[M]{a}}}
+\underbrace{F'\left({\sqrt[M]{a}}\right)\left(x-{\sqrt[M]{a}}\right)
+\ldots
+F\left({\sqrt[M]{a}}\right)^{\left(P\right)}
\frac{\left(x-{\sqrt[M]{a}}\right)^{P}}{P\,!}}_{=0}
\\
&+F\left(\xi\right)^{\left(P+1\right)}
\frac{\left(x-{\sqrt[M]{a}}\right)^{P+1}}{\left(P+1\right)!}\;,
\end{align*}
thus simplifies to
 \[
F\left(x\right)=
{\sqrt[M]{a}}
+F\left(\xi\right)^{\left(P+1\right)}
\frac{\left(x-{\sqrt[M]{a}}\right)^{P+1}}{\left(P+1\right)!}\;.
\]
Now plugging in $x=x_n$ sufficiently close to ${\sqrt[M]{a}}$ we obtain
\begin{align*}
F\left(x_n\right)=
{\sqrt[M]{a}}
+F\left(\xi_n\right)^{\left(P+1\right)}
\frac{\left(x_n-{\sqrt[M]{a}}\right)^{P+1}}{\left(P+1\right)!}
\end{align*}
Since $F\left(x_n\right)=x_{n+1}$ we get
\begin{align*}
x_{n+1}={\sqrt[M]{a}}
+F\left(\xi_n\right)^{\left(P+1\right)}
\frac{\left(x_n-{\sqrt[M]{a}}\right)^{P+1}}{\left(P+1\right)!}\,.
\end{align*}
Taking ${\sqrt[M]{a}}$ to the left side and dividing by
$\left(x_n-{\sqrt[M]{a}}\right)^{P+1}$ yields
\[
\frac{x_{n+1}-{\sqrt[M]{a}}}{\left(x_n-{\sqrt[M]{a}}\right)^{P+1}}=
F\left(\xi_n\right)^{\left(P+1\right)}
\frac{1}{\left(P+1\right)!}\;.
\]
As $x_n\rightarrow {\sqrt[M]{a}}$, since $\xi_n$ is trapped between $x_n$ and ${\sqrt[M]{a}}$ we conclude, by the continuity of $F\left(x\right)^{\left(P+1\right)}$ at ${\sqrt[M]{a}}$, that
\begin{align*}
\underset{n\to \infty}{\mathop{\lim }}\; \frac{x_{n+1}-{\sqrt[M]{a}}}{\left(x_n-{\sqrt[M]{a}}\right)^{P+1}}
=\underset{n\to \infty}{\mathop{\lim }}\;
F\left(\xi_n\right)^{\left(P+1\right)}
\frac{1}{\left(P+1\right)!}
=\frac{1}{\left(P+1\right)!}\,F\left({\sqrt[M]{a}}\right)^{\left(P+1\right)}
\end{align*}
We use an earlier result
\[
F{\left({\sqrt[M]{a}}\right)}^{\left({P+1}\right)}
=\frac{\left(-1\right)^P}{\sqrt[M]{a^P}}\prod_{\ell=1}^{P}\left(1+\ell\cdot M\right)\;.
\]
Plugging in gives
\begin{align*}
\underset{n\to \infty}{\mathop{\lim }}\;\;
\frac{x_{n+1}-{\sqrt[M]{a}}}{\left(x_n-{\sqrt[M]{a}}\right)^{P+1}}
=\frac{1}{\left(P+1\right)\,!}\cdot\frac{\left(-1\right)^P}{\sqrt[M]{a^P}}\prod_{\ell=1}^{P}\left(1+\ell\cdot M\right)\;.
\end{align*}
This shows that convergence is exactly of the order $\left(P+1\right)$, and
\[
\frac{1}{\left(P+1\right)\,!}\cdot\frac{\left(-1\right)^P}{\sqrt[M]{a^P}}\prod_{\ell=1}^{P}\left(1+\ell\cdot M\right)
\]
is the asymptotic error constant.
\section{Practical application of the method}
Let $a>0\,.$ We want to compute $\sqrt[M]{a}$. Additional we choose $P\in\mathbb{N}$. That gives order of convergence $\left(P+1\right)$.
\subsection{Polynomials up to order of convergence 5}
Before we begin with the practical application, let us take a look at the fixed point polynomials for $P \in \{1,2,3,4\}$. We obtain:

\begin{center}
\text{\underline{quadratic}}
\end{center}
\[
\left(1+\frac{1}{M}\right)\left( x-\frac{x^{M+1}}{a\cdot M}\right)
\]

\begin{center}
\text{\underline{cubic}}
\end{center}
\[
\left(1+\frac{1}{M}\right)\left(1+\frac{1}{2M}\right)\left(x-\frac{2\,x^{M+1}}{\left(M+1\right)a}+\frac{x^{2M+1}}{\left(2\,M+1\right)a^{2}}\right)
\]

\begin{center}
\text{\underline{quartic}}
\end{center}
\[
\left(1+\frac{1}{M}\right)\left(1+\frac{1}{2\,M}\right)\left(1+\frac{1}{3\,M}\right)\left(x-\frac{3\,x^{M+1}}{\left(M+1\right)a}+\frac{3\,x^{2\,M+1}}{\left(2\,M+1\right)a^{2}}-\frac{x^{3\,M+1}}{\left(3\,M+1\right)a^{3}}\right)
\]

\begin{center}
\text{\underline{quintic}}
\end{center}
\[
\left(1+\frac{1}{M}\right)\left(1+\frac{1}{2\,M}\right)\left(1+\frac{1}{3\,M}\right)\left(1+\frac{1}{4\,M}\right)\times
\]
\[
\left(x-\frac{4\,x^{M+1}}{\left(M+1\right)a}+\frac{6\,x^{2\,M+1}}{\left(2\,M+1\right)a^{2}}-\frac{4\,x^{3\,M+1}}{\left(3\,M+1\right)a^{3}}+\frac{x^{4\,M+1}}{\left(4\,M+1\right)a^{4}}\right)
\]

\subsection{General description of the method}
Let $a>0$. $M\in\mathbb{N}$ and  $P\in\mathbb{N}$. We start with the function
\[f\left(t\right)=1-\frac{t^M}{a}\;.\]
First we evaluate the Integral in order to get the fixed point function
\[
F\left(x\right)
=\prod_{\ell=1}^{P}\left(1+\frac{1}{\ell\cdot M}\right)
\int\limits_{0}^{x}\!f\left(t\right)^{P}{\rm d}t
=\prod_{\ell=1}^{P}\left(1+\frac{1}{\ell\cdot M}\right)
\int\limits_{0}^{x}\!\left(1-{\frac{{t}^{M}}{a}}\right)^{P}{\rm d}t\;.
\]
Applying the binomial theorem to
\[\left(1-\frac{t^M}{a}\right)^{P}\]
gives
\begin{align*}
\left(1-\frac{t^M}{a}\right)^{P}
&=\sum\limits_{k=0}^{P}\left(-1\right)^k\,\cdot\binom{P}{k}\,\left(\frac{t^M}{a}\right)^{\,k}
\\
&=\sum\limits_{k=0}^{P}\left(-1\right)^k\,\cdot\binom{P}{k}\cdot\frac{t^{\,k\,\cdot M}}{a^{\,k}}
\\
&=\sum\limits_{k=0}^{P}\frac{\left(-1\right)^{\,k}}{a^{\,k}}\cdot\binom{P}{k}\cdot t^{\,k\,\cdot M}\;.
\end{align*}
This sum can be easily integrated:
\begin{align*}
\int\limits_{0}^{x}\!\left(1-{\frac{{t}^{M}}{a}}\right)^{P}{\rm d}t
&=\int\limits_{0}^{x}\!\sum\limits_{k=0}^{P}\frac{\left(-1\right)^{\,k}}{a^{\,k}}\cdot\binom{P}{k}\cdot x^{\,k\,\cdot M}{\rm d}t
\\
&=\sum\limits_{k=0}^{P}\frac{\left(-1\right)^{\,k}}{a^{\,k}}\cdot\binom{P}{k}\cdot\frac{x^{\,k\,\cdot M+1}}{k\,\cdot M+1}\;.
\end{align*}
Now we just need to multiply by the product
\[
\prod_{\ell=1}^{P}\left(1+\frac{1}{\ell\cdot M}\right)
\]
and the fixed point function is complete:
\[
F\left(x\right)=\prod_{\ell=1}^{P}\left(1+\frac{1}{\ell\cdot M}\right)
\sum\limits_{k=0}^{P}\frac{\left(-1\right)^{\,k}}{a^{\,k}}\cdot\binom{P}{k}\cdot\frac{x^{\,k\,\cdot M+1}}{k\,\cdot M+1}\;.
\]
This is a polynomial of degree $\left(P\cdot M+1\right)$ with $\left(P+1\right)$ terms.
$ \\ $
$ \\ $
We observe that the somewhat complicated coefficients
\[
\prod_{\ell=1}^{P}\left(1+\frac{1}{\ell\cdot M}\right)
\frac{\left(-1\right)^{\,k}}{a^{\,k}}\cdot\binom{P}{k}\cdot\frac{1}{k\,\cdot M+1}\;;\;k=0,1,\ldots,P\]
do not change throughout the calculation. It is therefore advantageous to calculate them beforehand and store them in a vector:
\[
c_{\,k}=\prod_{\ell=1}^{P}\left(1+\frac{1}{\ell\cdot M}\right)
\frac{\left(-1\right)^{\,k}}{a^{\,k}}\cdot\binom{P}{k}\cdot\frac{1}{k\,\cdot M+1}\;;\;k=0,1,\ldots,P\;.\]
The fixed point function is then simplified to
\[F\left(x\right)=\sum\limits_{k=0}^{P}c_{\,k}\cdot x^{\,k\,\cdot M+1}\]
\underline{Note:} The calculation of $F\left(x\right)$ is a polynomial evaluation. There are efficient methods for polynomial computing (Horner's scheme).
$ \\ $
$ \\ $
For practical computations, it is useful to shift the index. Instead of
\[x_{n+1}=\sum\limits_{k=0}^{P}c_{\,k}\cdot x_{n}^{\,k\,\cdot M+1}\;;\;n=0,1,\ldots\]
we will write
\[x_{n}=\sum\limits_{k=0}^{P}c_{\,k}\cdot x_{n-1}^{\,k\,\cdot M+1}\;;\;n=1,2,\ldots\]
$ \\ $
Let $x_0\in\mathbb{R}$ an initial value sufficiently close to$\sqrt[M]{a}$.
$ \\ $
Now we can start the iteration procedure

\[\underline{step\;\;1}\]
\[x_1=\sum\limits_{k=0}^{P}c_{\,k}\cdot x_0^{\,k\,\cdot M+1}\;;\]

\[\underline{step\;\;2}\]
\[x_2=\sum\limits_{k=0}^{P}c_{\,k}\cdot x_1^{\,k\,\cdot M+1}\;;\]

\[\underline{step\;\;3}\]
\[x_3=\sum\limits_{k=0}^{P}c_{\,k}\cdot x_2^{\,k\,\cdot M+1}\;;\]
%
%
%
\[\vdots\]
\[\underline{Step\;\;n}\]
\[x_{n}=\sum\limits_{k=0}^{P}c_{\,k}\cdot x_{n-1}^{\,k\,\cdot M+1}\;;\]
$ \\ $
If, now $n\rightarrow \infty$, then $x_{n}\rightarrow \sqrt[M]{a}$ with order of convergence $\left(P+1\right)$. With unlimited computing time and unlimited storage space, we can calculate $\sqrt[M]{a}$  with arbitrary precision. Unfortunately, we do not have anything like that available. We must be content with a finite number of iterations and a finite value for $\sqrt[M]{a}$.
$ \\ $
We require a termination criterion. For this purpose, we specify $\epsilon>0$ and perform the iteration until the following applies for two consecutive values for the first time
\[\left|x_n-x_{n-1}\right|<\epsilon\;.\]
\subsection{Computational results}
\textbf{{Example $1$\,: Computation of $\sqrt[3]{10}$}}
$ \\ $
In the first example, we calculate $\sqrt[3]{10}$. So we have $a=10$ and $M=3$.  Plugging in $a=10\;,M=3\;$ yields $f\left(t\right)=1-\frac{t^3}{10}\;.$
We choose $P=1$. This gives the order of convergence $\left(P+1\right)=2$ (quadratically). Each additional Iteration approximately doubles the number of correct digits. The fixed point polynomial is then given by evaluating the integral
\begin{align*}
F\left(x\right)
=\prod_{k=1}^{1}\left(1+\frac{1}{3\cdot k}\right)
\int\limits_{0}^{x}\!\underbrace{\left(1-{\frac{{t}^{3}}{10}}\right)}_{f\left(t\right)}{dt}
=\left(1+\frac{1}{3}\right)\left(x-\frac{1}{40}\,x^4\right)
=\frac{4}{3}\,x-\frac{1}{30}\,x^4\;.\end{align*}
Let's look at the functions
together with $y=x\,:$
\begin{center}
\includegraphics[width=0.95\linewidth]{"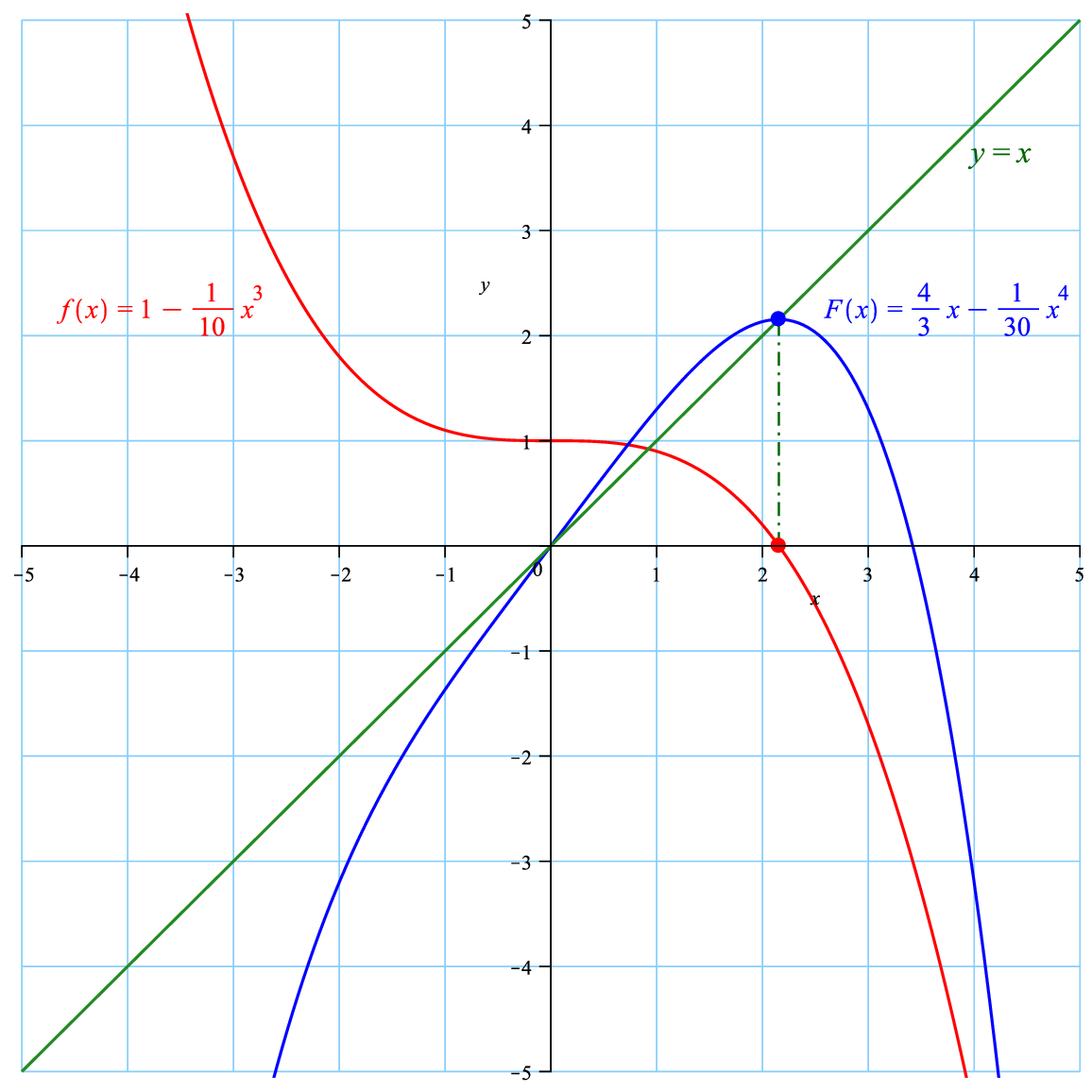"}
\textbf{\underline{figure 1}}
\end{center}
With initial value $x_0=2$, we compute the first six iterations

\[\underline{step\;\;1}\]
\[x_{1}=\frac{4}{3}x_{0}-\frac{1}{30}x_{0}^{4}\]
\[x_{1}=2.133333333333333333333333333333333333333\]
\[{|x_{1}-x_{0}|}=0.1333333333333333333333333333333333333333\]
\[\underline{step\;\;2}\]
\[x_{2}=\frac{4}{3}x_{1}-\frac{1}{30}x_{1}^{4}\]
\[x_{2}=2.154024032921810699588477366255144032922\]
\[{|x_{2}-x_{1}|}=0.02069069958847736625514403292181069958848\]
\[\underline{step\;\;3}\]
\[x_{3}=\frac{4}{3}x_{2}-\frac{1}{30}x_{2}^{4}\]
\[x_{3}=2.154434533500953092649669501763572523986\]
\[{|x_{3}-x_{2}|}=0.0004105005791423930611921355084284910642133\]
\[\underline{step\;\;4}\]
\[x_{4}=\frac{4}{3}x_{3}-\frac{1}{30}x_{3}^{4}\]
\[x_{4}=2.154434690031860976181374509716973801410\]
\[{|x_{4}-x_{3}|}={1.565309078835317050079534012774237318926\times10^{-7}}\]
\[\underline{step\;\;5}\]
\[x_{5}=\frac{4}{3}x_{4}-\frac{1}{30}x_{4}^{4}\]
\[x_{5}=2.154434690031883721759293566039074794849\]
\[{|x_{5}-x_{4}|}={2.274557791905632210099343907978738060749\times10^{-14}}\]
\[\underline{step\;\;6}\]
\[x_{6}=\frac{4}{3}x_{5}-\frac{1}{30}x_{5}^{4}\]
\[x_{6}=2.154434690031883721759293566519350495259\]
\[{|x_{6}-x_{5}|}={4.802757004105093077094334087308664908888\times10^{-28}}\]
$ \\ $
The differences $\left|x_5-x_4\right|$ and $\left|x_6-x_5\right|$ clearly show the order of convergence $2$. The negative exponents double with each iteration step. 
$ \\ $
$ \\ $
$ \\ $
\textbf{\underline{Example $2$\,: Computing one million digits of $\sqrt{2}$}}
$ \\ $
The next example is a little more challenging. We calculate $\sqrt{2}$ to at least one million digits. That means: $a=2$, $M=2$ and $f\left(t\right)=1-\frac{t^2}{2}\,.$ We choose $P=3$ and thus obtain order of convergence $3+1=4$. As the starting value, we use $x_0=1.414213562373$, accurate to 16 digits. Since the number of digits quadruples with each iteration step, we hope to reach the desired result after nine iterations. We begin with constructing the fixed point function
\[
F\left(x\right)
=\prod_{\ell=1}^{3}\left(1+{\frac{1}{2\,\ell}}\right)\int\limits_{0}^{x}\!\left(1-{\frac{{x}^{2}}{2}}\right)^{3}\,{\rm d}t
=\ldots=\frac{35}{16}\,x-\frac{35}{32}\,{x}^{3}+\frac{21}{64}\,{x}^{5}-\frac{5}{128}\,{x}^{7}\;.
\]
Let's look at the functions together with $y=x$.
\begin{center}
\includegraphics[width=0.95\linewidth]{"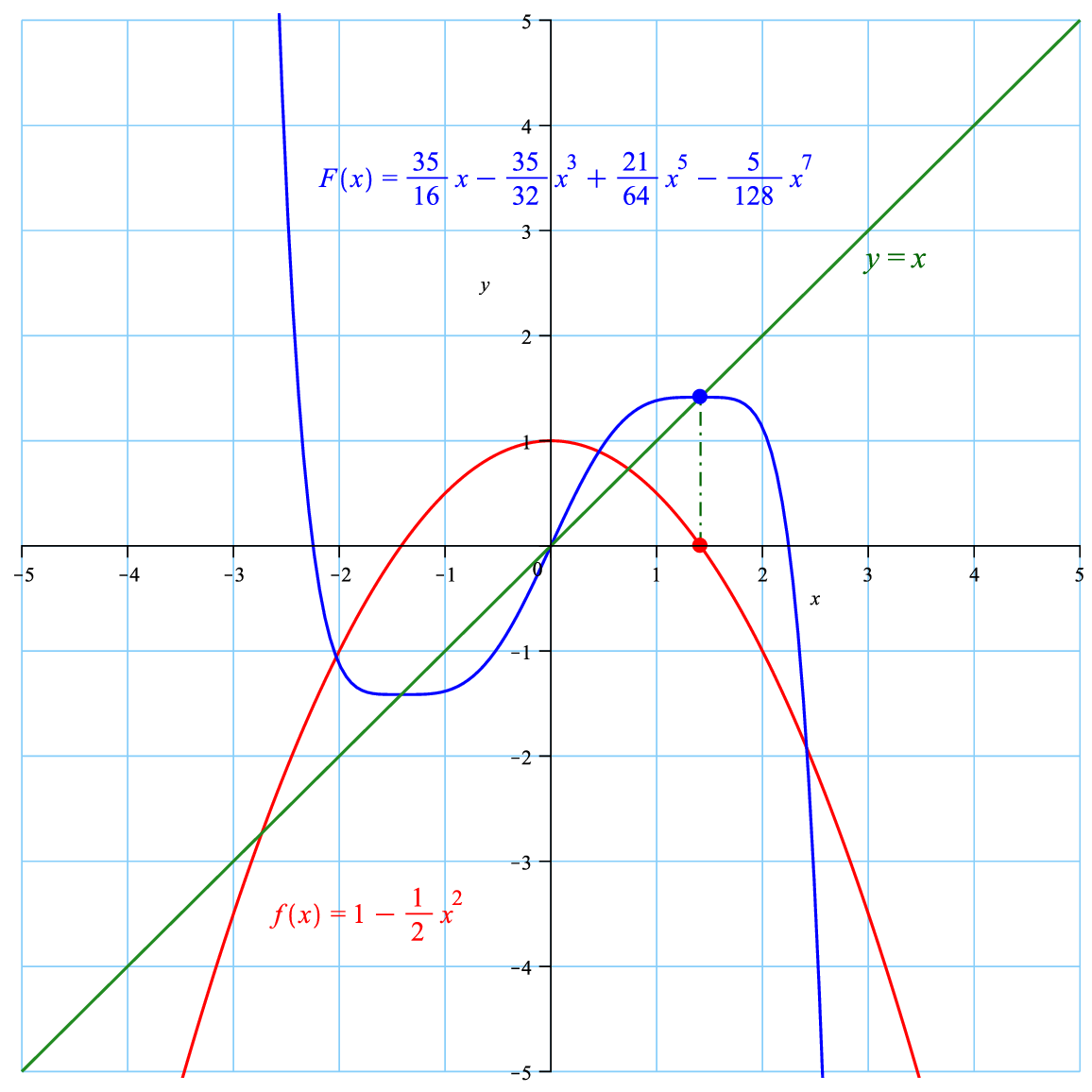"}
\end{center}
\begin{center}
\textbf{\underline{figure 2}}
\end{center}
$ \\ $
With the initial value $x_0=1.414213562373$, we obtain the sequence
\[\underline{step\;\;1}\]
\[x_{1}=\frac{35}{16}x_{0}-\frac{35}{32}x_{0}^{3}+\frac{21}{64}x_{0}^{5}-\frac{5}{128}x_{0}^{7}\]
\[x_{1}=1.414213562373095048801688724209698078570\]
\[{|x_{1}-x_{0}|}={4.880168872420969807856967187537694807318\times10^{-17}}\]
\[\underline{step\;\;2}\]
\[x_{2}=\frac{35}{16}x_{1}-\frac{35}{32}x_{1}^{3}+\frac{21}{64}x_{1}^{5}-\frac{5}{128}x_{1}^{7}\]
\[x_{2}=1.414213562373095048801688724209698078570\]
\[{|x_{2}-x_{1}|}={8.773491625654111352087407579690431191435\times10^{-66}}\]
\[\underline{step\;\;3}\]
\[x_{3}=\frac{35}{16}x_{2}-\frac{35}{32}x_{2}^{3}+\frac{21}{64}x_{2}^{5}-\frac{5}{128}x_{2}^{7}\]
\[x_{3}=1.414213562373095048801688724209698078570\]
\[{|x_{3}-x_{2}|}={9.164798637556653681657805406878049888878\times10^{-261}}\]
\[\underline{step\;\;4}\]
\[x_{4}=\frac{35}{16}x_{3}-\frac{35}{32}x_{3}^{3}+\frac{21}{64}x_{3}^{5}-\frac{5}{128}x_{3}^{7}\]
\[x_{4}=1.414213562373095048801688724209698078570\]
\[{|x_{4}-x_{3}|}={1.091251298365935101705686744387078883102\times10^{-1040}}\]
\[\underline{step\;\;5}\]
\[x_{5}=\frac{35}{16}x_{4}-\frac{35}{32}x_{4}^{3}+\frac{21}{64}x_{4}^{5}-\frac{5}{128}x_{4}^{7}\]
\[x_{5}=1.414213562373095048801688724209698078570\]
\[{|x_{5}-x_{4}|}={2.193472316487722705810599621121648551289\times10^{-4160}}\]
\[\underline{step\;\;6}\]
\[x_{6}=\frac{35}{16}x_{5}-\frac{35}{32}x_{5}^{3}+\frac{21}{64}x_{5}^{5}-\frac{5}{128}x_{5}^{7}\]
\[x_{6}=1.414213562373095048801688724209698078570\]
\[{|x_{6}-x_{5}|}={3.580648536099876136173035995717511426715\times10^{-16639}}\]
\[\underline{step\;\;7}\]
\[x_{7}=\frac{35}{16}x_{6}-\frac{35}{32}x_{6}^{3}+\frac{21}{64}x_{6}^{5}-\frac{5}{128}x_{6}^{7}\]
\[x_{7}=1.414213562373095048801688724209698078570\]
\[{|x_{7}-x_{6}|}={2.542610528450840832485991523758935060375\times10^{-66554}}\]
\[\underline{step\;\;8}\]
\[x_{8}=\frac{35}{16}x_{7}-\frac{35}{32}x_{7}^{3}+\frac{21}{64}x_{7}^{5}-\frac{5}{128}x_{7}^{7}\]
\[x_{8}=1.414213562373095048801688724209698078570\]
\[{|x_{8}-x_{7}|}={6.464760315447686077979797373449536529093\times10^{-266215}}\]
\[\underline{step\;\;9}\]
\[x_{9}=\frac{35}{16}x_{8}-\frac{35}{32}x_{8}^{3}+\frac{21}{64}x_{8}^{5}-\frac{5}{128}x_{8}^{7}\]
\[x_{9}=1.414213562373095048801688724209698078570\]
\[{|x_{9}-x_{8}|}={2.701735162639912537134047073288055961734\times10^{-1064857}}\]
$\\ $
The result is available within 3 seconds. With nine iterations, we actually obtained over a million digits for $\sqrt{2}$. From the negative exponents of the differences, we can see the convergence order 4. The negative exponent approximately quadruples with each step.
$\\ $
$\\ $
We performed all computational results by using a homemade PC with the following hardware configuration: Motherboard ASUS PRIME A320M-K with CPU AMD Ryzen 5 5600G 6 CORE 3.90-4.40 GHz and 32 GB RAM. The used software was MAPLE 2025.2 by  Maplesoft, Waterloo Maple Inc. This software was also used to assist with formatting the mathematical equations into LaTeX.

\newpage
\selectlanguage{ngerman}
\begin{center}
\LARGE{\textbf{Ein iteratives Wurzelberechnungsverfahren \\ mit beliebiger Konvergenzordnung}}
$$
$$
\large
{Alois Schiessl}
\footnote[1]{University of Regensburg\;;\;E-Mail: \texttt{aloisschiessl@web.de}}
\vskip 0.75cm
\centerline{\today}
\end{center}
\[\]
\[\]
\begin{abstract}
Es sei $a\in \mathbb{R}^{+}\backslash\left\{0\right\}$ und $M\in\mathbb{N}$. Vorgelegt ist die Gleichung $1-\frac{t^M}{a}=0$. Die reelle Lösung hiervon ist $\sqrt[M]{a}$. In dieser Veröffentlichung stellen wir ein Verfahren vor, das die Berechnung von $\sqrt[M]{a}$ mit beliebiger Konvergenzordnung $P\in\mathbb{N}$ ermöglicht. Wir definieren die Fixpunktfunktion
\[F\left(x\right)
=\prod_{\ell=1}^{P}\left(1+\frac{1}{\ell\cdot M}\right)
\int\limits_{0}^{x}\!\left(1-{\frac{{t}^{M}}{a}}\right)^{P}{\rm d}t
=\sum\limits_{k=0}^{P}\frac{\left(-1\right)^{\,k}}{a^{\,k}}\cdot\binom{P}{k}\cdot\frac{x^{\,k\,\cdot M+1}}{k\,\cdot M+1}\;.
\]
Das ist ein Polynom vom Grad $\left(P\cdot M+1\right)$ mit $\left(P+1\right)$ Summanden. Die Berechnung von $\sqrt[M]{a}$ reduziert sich dadurch zu einer Polynomauswertung. Hierfür gibt es effiziente Verfahren.
\end{abstract}
\setcounter{section}{0}
\section{Einleitung und Hauptergebnis}
Die wohl bekannteste Methode um zu einer positiven reellen Zahl $a$ und einer natürlichen Zahl $M$ einen Näherungswert für $\sqrt[M]a$ zu berechnen besteht darin, auf die Funktion $f\left(x\right)=x^M-a$ das Newton-Verfahren
\[
x_{n+1}=x_n-\frac{f\left(x_n\right)}{f\,'\left(x_n\right)}
\]
anzuwenden. Das führt zur Fixpunkt-Iteration
\[x_{n+1}
=x_n-\frac{x_n^{M}-a}{M\,x_n^{M -1}}
=x_n-\frac{1}{M}\left(x_n-\frac{a}{x_n^{M -1}}\right)\;;\;n=0,1,\ldots\]
Liegt der Startwert $x_0$ hinreichend nahe bei $\sqrt[M]a$ so konvergiert die Folge gegen $\sqrt[M]a$ mit Konvergenzordnung $2$ (quadratische Konvergenz). Das bedeutet, dass sich die Anzahl der gültigen Dezimalstellen mit jedem Schritt näherungsweise verdoppelt.
$ \\ $
Wir geben hier ein Verfahren an, das nur mit Polynomen arbeitet und beliebige Konvergenzordnung ermöglicht. Das Verfahren basiert auf der reellen Funktion
\[f\left(t\right)=1-\frac{t^M}{a}\;.\]
Wir formulieren das folgende Theorem
\begin{theorem}
Es sei $a\in \mathbb{R}\backslash\left\{0\right\}$ und $M\in \mathbb{N}$. Wir definieren die Funktion
\[f: \mathbb{R} \to \mathbb{R}\]
\[t\mapsto f\left(t\right)\]
\[f\left(t\right)=1-\frac{t^M}{a}\]
$ \\ $
Weiter definieren wir die Integral-Funktion
\[F: \mathbb{R} \to \mathbb{R}\]
\[x\mapsto F\left(x\right)\]
\[F\left(x\right)
=\prod_{\ell=1}^{P}\left(1+\frac{1}{\ell\cdot M}\right)\int\limits_{0}^{x}\!f\left(t\right)^{P}{\rm d}t
=\prod_{\ell=1}^{P}\left(1+\frac{1}{\ell\cdot M}\right)
\int\limits_{0}^{x}\!\left(1-{\frac{{t}^{M}}{a}}\right)^{P}{\rm d}t
\]
Es sei $a>0.$ Dann gelten für $F\left(x\right)$ die folgenden Aussagen
\begin{itemize}
\item[(a)]
$\sqrt[M]{a}\,$ ist Fixpunkt, d.h. es gilt: $F\left(\sqrt[M]{a}\right)=\sqrt[M]{a}\,.$
\item[(b)]
Die Ableitungen im Fixpunkt $\sqrt[M]{a}$ sind gegeben durch
\begin{align*}
F{\left(\sqrt[M]{a}\right)}^{\left(k\right)}=\left\{
\begin{aligned}
0\qquad\qquad\qquad\qquad\qquad\;\; 1\,\leq &\,k\,\leq {P}
\\
\\
\frac{(-1)^P}{\sqrt[M]{a^P}}\prod_{\ell=1}^{P}\left(1+\ell\cdot M\right)
\qquad\qquad\quad k=&P+1
\\ 
\end{aligned} 
\right.
\end{align*}
\item[(c)]
Für jeden Startwert $x_0\in\mathbb{R}$ hinreichend nahe bei $\sqrt[M]{a}$ konvergiert die Folge
\[
x_{n+1}=F\left(x_{{n}}\right)
\]
gegen $\sqrt[M]{a}$ mit Konvergenzordnung genau $\left(P+1\right)$.
\end{itemize}
\end{theorem}
\section{Beweis des Theorem}
In diesem Kapitel weisen wir die Aussagen zum Theorem nach. Die Beweise sind nicht schwierig. Es werden nur Lehrsätze der klassichen Analysis benötigt.
\subsection{Fixpunkt}
Wir beginnen mit der Fixpunktaussage. Wir haben zu beweisen, dass für $a>0$ die folgende Aussage gilt
\[F\left(\sqrt[M]{a}\right)=\prod_{\ell=1}^{P}\left(1+\frac{1}{\ell\cdot M}\right)
\int\limits_{0}^{\sqrt[M]{a}}\!\left(1-{\frac{{t}^{M}}{a}}\right)^{P}{\rm d}t=\sqrt[M]{a}\,.\]
Wir befassen uns zunächst mit dem bestimmten Integral
\[
\int\limits_{0}^{\sqrt[M]{a}}\!\left(1-\frac{{t}^{M}}{a}\right)^{P}{\rm d}t\,.\]
Als erstes substituieren wir
$\frac{t^M}{a}=u$ und lösen nach $t$ auf. Das ergibt $t=\sqrt[M]{a\cdot u}=\sqrt[M]{a}\cdot\sqrt[M]{u}$. Wir erhalten neue Integrationsgrenzen:
\\
$t=0$ ergibt ebenfalls $u=0$ und für $t=\sqrt[M]{a}$ erhalten wir $u=1$. Wir benötigen noch die Variablentransformation
\[\frac{{\rm d} \,t} {{\rm d}\, u}=\sqrt[M]{a}\cdot\frac{1}{M}\cdot u^{\frac{1}{M}-1}\,.\]
Wir erhalten damit ein neues Integral
\[
\int\limits_{0}^{\sqrt[M]{a}}\!\left(1-{\frac{{t}^{M}}{a}}\right)^{P}{\rm d}t
=\int\limits_{0}^{1}\!\left(1-u\right)^{P}\sqrt[M]{a}\cdot\frac{1}{M}\cdot u^{\frac{1}{M}-1}{\rm d}u
=\frac{\sqrt[M]{a}}{M}\int\limits_{0}^{1}\!u^{\frac{1}{M}-1}\cdot\left(1-u\right)^{P}{\rm d}u\,.
\]
Die nächste Aufgabe besteht jetzt darin das Integral auf der rechten Seite auszuwerten. Wir erinnern uns an die reelle Beta-Funktion \cite{AS}
\[{\rm B}\left(\alpha,\beta\right)=\int\limits_{0}^{1}\!{u}^{\alpha-1}\left(1-u\right)^{\beta-1}\,{\rm d}u
\;;\,\alpha>0\;;\,\beta>0\,.\]
Diese ist nahe verwandt mit der Gammafunktion \cite{EA}
\[{\rm B}\left(\alpha,\beta\right)
={\frac{\Gamma\left(\alpha\right)\Gamma\left(\beta\right)}{\Gamma\left(\alpha+\beta\right)}}\,.\]
Wir substituieren
\[\alpha=\frac{1}{M}\,;\]
\[\beta=P+1\]
und erhalten eine Darstellung durch die $\Gamma$-Funktion
\[
\int\limits_{0}^{1}\!u^{\frac{1}{M}-1}\cdot\left(1-u\right)^{P}{\rm d}u
=\frac{\Gamma\left(\frac{1}{M}\right)\cdot\Gamma\left(P+1\right)}{\Gamma\left(P+1+\frac{1}{M}\right)}\,.
\]
Nun ist $P$ eine natürliche Zahl, somit gilt $\Gamma\left(P+1\right)=P\,!$ und der Term  vereinfacht sich zu
\[
\frac{\Gamma\left(\frac{1}{M}\right)\cdot\Gamma\left(P+1\right)}{\Gamma\left(P+1+\frac{1}{M}\right)}
=\frac{\Gamma\left(\frac{1}{M}\right)P\,!}{\Gamma\left(P+1+\frac{1}{M}\right)}\,.
\]
Als nächstes befassen wir uns mit
\[
\Gamma\left(P+1+\frac{1}{M}\right)\,.
\]
Durch fortgesetzte Anwendung der Definition der Gammafunktion
\[\Gamma\!\left(1+x\right)=x\,\Gamma\!\left(x\right)\,;\]
 mit $x=\frac{1}{M}$ erhalten wir nacheinander
\[\Gamma\!\left(1+\frac{1}{M}\right)=\frac{\Gamma\!\left(\frac{1}{M}\right)}{M}\,;\]
\begin{align*}
\Gamma\!\left(2+\frac{1}{M}\right)
=\prod\limits_{\ell=1}^{1}\left(1+\frac{1}{\ell\cdot M}\right)\frac{1\,!\cdot\Gamma\left(\frac{1}{M}\right)}{M} \,;
\end{align*}
\begin{align*}
\Gamma\!\left(3+\frac{1}{M}\right)
=\prod\limits_{\ell=1}^{2}\left(1+\frac{1}{\ell\cdot M}\right)\frac{2\,!\cdot\Gamma\left(\frac{1}{M}\right)}{M} \,;
\end{align*}
\begin{align*}
\Gamma\!\left(4+\frac{1}{M}\right)
=\prod\limits_{\ell=1}^{3}\left(1+\frac{1}{\ell\cdot M}\right)\frac{3\,!\cdot\Gamma\left(\frac{1}{M}\right)}{M}\,;
\end{align*}
\[\vdots\]
Es ist eine einfache Angelegenheit durch vollständiger Induktion folgende Aussage zu beweisen:
\[
\Gamma\left(P+1+\frac{1}{M}\right)
=\prod\limits_{\ell=1}^{P}\left(1+\frac{1}{\ell\cdot M}\right)\frac{P\,!\cdot\Gamma\left(\frac{1}{M}\right)}{M}\,.
\]
Für das bestimmte Integral erhalten wir
\[
\int\limits_{0}^{1}\!u^{\frac{1}{M}-1}\cdot\left(1-u\right)^{P}{\rm d}u
=\frac{\Gamma\left(\frac{1}{M}\right)\cdot\Gamma\left(P+1\right)}{\Gamma\left(P+1+\frac{1}{M}\right)}
=\frac{\Gamma\left(\frac{1}{M}\right)\cdot P\,!}{\prod\limits_{\ell=1}^{P}\left(1+\frac{1}{\ell\cdot M}\right)\frac{P\,!\cdot\Gamma\left(\frac{1}{M}\right)}{M}}
=\frac{M}{\prod\limits_{\ell=1}^{P}\left(1+\frac{1}{\ell\cdot M}\right)}\,.
\]
Jetzt brauchen wir nur noch alles einzusetzen. Das ergibt dann
\begin{align*}
F\left(\sqrt[M]{a}\right)
&=\prod_{\ell=1}^{P}\left(1+\frac{1}{\ell\cdot M}\right)\int\limits_{0}^{\sqrt[M]{a}}\!\left(1-{\frac{{t}^{M}}{a}}\right)^{P}{\rm d}t
\\
&=\prod_{\ell=1}^{P}\left(1+\frac{1}{\ell\cdot M}\right)\frac{\sqrt[M]{a}}{M}\int\limits_{0}^{1}\!u^{\frac{1}{M}-1}\cdot\left(1-u\right)^{P}{\rm d}u
\\
&=\prod_{\ell=1}^{P}\left(1+\frac{1}{\ell\cdot M}\right)\frac{\sqrt[M]{a}}{M}\cdot \frac{M}{\prod\limits_{\ell=1}^{P}\left(1+\frac{1}{\ell\cdot M}\right)}
\\
&=\prod_{\ell=1}^{P}\left(1+\frac{1}{\ell\cdot M}\right)\frac{\sqrt[M]{a}}{\prod\limits_{\ell=1}^{P}\left(1+\frac{1}{\ell\cdot M}\right)}
=\sqrt[M]{a}\,.
\end{align*}
Damit ist die Fixpunkteigenschaft bewiesen.
\subsection{Ableitungen}
Als nächstes wenden wir uns den Ableitungen zu.
\subsubsection{erste Ableitung}
Wir beginnen mit der ersten Ableitung. Diese lässt sich wegen der speziellen Bauweise von $F\left(x\right)$ besonders einfach berechnen. Durch Differentiation erhalten wir sofort
\[
F\,'\left(x\right)=\prod_{\ell=1}^{P}\left(1+\frac{1}{\ell\cdot M}\right)\left(1-{\frac{{x}^{M}}{a}}\right)^{P}\,.
\]
Eingesetzt $x=\sqrt[M]{a}$ ergibt $F\,'\left(\sqrt[M]{a}\right)=0$. Dieses Ergebnis benötigen wir später.
\subsubsection{höhere Ableitungen}
Wir haben zu beweisen, dass für die Ableitungen bis $\left(P+1\right)$ im Fixpunkt $\sqrt[M]{a}$ folgendes gilt:
\begin{align*}
F{\left(\sqrt[M]{a}\right)}^{\left(k\right)}=\left\{
\begin{aligned}
0\qquad\qquad\qquad\qquad 1\,\le &\,k\,\leq{P}
\\
\\
\frac{\left(-1\right)^P}{\sqrt[M]{a^P}}\prod_{\ell=1}^{P}\left(1+\ell\cdot M\right)
\qquad\qquad\quad k=&P+1
\\ 
\end{aligned} 
\right.
\end{align*}
Dazu müssen wir $F\left(x\right)$ mal $\left(P+1\right)$ differenzieren und dann $x=\sqrt[M]{a}$ einsetzen. Da die erste Ableitung bereits vorliegt vereinfacht sich die Aufgabe:$F\,'\left(x\right)$ 
\begin{align*}
\left.F\left(x\right)^{\left(P+1\right)}\right|_{x={\sqrt[M]{a}}}
=\left.\bigl(F\,'\left(x\right)\bigr)^{\left(P\right)}\right|_{x={\sqrt[M]{a}}}
&=\left.\left(\prod_{k=1}^{P}\left(1+{\frac{1}{k\,\cdot M}}\right)\cdot\left(1-\frac{x^M}{a}\right)^{P}\right)^{\left(P\,\right)}\right|_{x={\sqrt[M]{a}}}
\\
&=\prod_{k=1}^{P}\left(1+{\frac{1}{k\,\cdot M}}\right)\cdot\left.\left(\left(1-\frac{x^M}{a}\right)^{P}\right)^{\left(P\,\right)}\right|_{x={\sqrt[M]{a}}}\,.
\end{align*}
Auf der rechten Seite benötigen wir im Fixpunkt $\sqrt[M]{a}$ die $P^{\,te}$ Ableitung von
\[\left(1-\frac{x^M}{a}\right)^P\;.\]
Es bezeichne
\[f\left(x\right)=\left(1-\frac{x^M}{a}\right).\]
Vorab berechnen wir einige spezielle Werte, die wir später benötigen:
\[f\left(\sqrt[M]{a}\right)=0\;;\]
\[f'\left(\sqrt[M]{a}\right)=\ldots=\frac{-M}{\sqrt[M]{a}}\;;\]
Als nächstes entwickeln wir $f\left(x\right)$ in eine Taylorreihe erster Ordnung um $\sqrt[M]{a}:$ 
\[
f\left(x\right)
=f\left(\sqrt[M]{a}\right)
+f'\left(\sqrt[M]{a}\right)\left(x-\sqrt[M]{a}\right)
+\frac{f''\left(\xi\right)}{2}\left(x-\sqrt[M]{a}\right)^{2},
\]
wobei $\xi$ zwischen $\sqrt[M]{a}$ and $x$ liegt. Wir setzen die oben berechneten speziellen Werte ein und erhalten die einfache Darstellung:
\[
f\left(x\right)
=
\frac{-M}{\sqrt[M]{a}}\left(x-\sqrt[M]{a}\right)
+\frac{f''\left(\xi\right)}{2}\left(x-\sqrt[M]{a}\right)^{2}\,.
\]
Wir können $\left(x-\sqrt[M]{a}\right)$ ausklammern
\[
f\left(x\right)=
\left(x-\sqrt[M]{a}\right)
\left(\frac{-M}{\sqrt[M]{a}}+\frac{f''\left(\xi\right)}{2}\left(x-\sqrt[M]{a}\right)\right)\,.
\]
Jetzt erheben wir $f\left(x\right)$ in die $P^{\,te}$ Potenz
\[
f\left(x\right)^{P}=
\left(x-\sqrt[M]{a}\right)^P
\left(\frac{-M}{\sqrt[M]{a}}+\frac{f''\left(\xi\right)}{2}\left(x-\sqrt[M]{a}\right)\right)^P.
\]
Die höheren Ableitungen berechnen wir mit der allgemeinen Leibniz-Regel \cite{HJ}, die besagt:
\\
Seien  $u\left(x\right)$ und $v\left(x\right)$ zwei über demselben Intervall $n$-mal differenzierbare Funktionen, dann ist das Produkt $u\left(x\right)\cdot v\left(x\right)$ ebenfalls $n$-mal differenzierbar und es gilt: 
\[
\bigl(u\left(x\right)\cdot v\left(x\right)\bigr)^{\left(n\right)}
=\sum\limits_{k=0}^{n}\binom{n}{k}\,u\left(x\right)^{\left(k\right)}\cdot v\left(x\right)^{\left(n-k\right)}.
\]
Hierbei wollen wir unter $u\left(x\right)^{\left(0\right)}=u\left(x\right)$ und $v\left(x\right)^{\left(0\right)}=v\left(x\right)$ verstehen.
$\\ $
$\\ $
In unserem Fall wählen wir:
\[
n=P\;;
\]
\[
u\left(x\right)=\left(x-\sqrt[M]{a}\right)^{P}\;;
\]
\[
v\left(x\right)=\left(\frac{-M}{\sqrt[M]{a}}+\frac{f''\left(x\right)}{2}\left(x-\sqrt[M]{a}\right)\right)^{P},
\]
und erhalten
\[
f\left(x\right)^{P}=
\underbrace{{{\left(x-\sqrt[M]{a}\right)}^{P}}}_{u\left(x\right)}\cdot
\underbrace{{{\left(\frac{-M}{\sqrt[M]{a}}+\frac{f''{{\left(\xi\right)}}}{2}\left(x-\sqrt[M]{a}\right)\right)}^{P}}}_{v\left(x\right)}.
\]
Hierauf wenden wir die Leibniz-Regel an
\[
\left(f\left(x\right)^{P}\right)^{\left({P}\right)}
=\sum\limits_{k=0}^{{P}}\binom{{P}}{k}
\,\left(\left(x-\sqrt[M]{a}\right)^{P}\right)^{\left(k\right)}
\times
\left(\left(\frac{-M}{\sqrt[M]{a}}+\frac{f''{\left(\xi\right)}}{2}\left(x-\sqrt[M]{a}\right)\right)^{{P}}\right)^{\left({P-k}\right)}.
\]
$ \\ $
Wir sehen uns zuerst die Folge der Ableitungen an
\[
u\left(x\right)^{\left(k\right)}
=\Bigl(\left(x-\sqrt[M]{a}\right)^{P}\Bigr)^{\left(k\right)},\;k=0,1,\ldots,P\;.
\]
Aus einer Formelsammlung \cite{formula} verwenden wir die folgende Aussage
\begin{align*}
\Bigl(\left(x-\sqrt[M]{a}\right)^{P}\Bigr)^{\left(k\right)}
=\frac{P\,!}{\left(P-k\right)\,!}\left(x-\sqrt[M]{a}\right)^{P-k}
,\;k=0,1,\ldots,P\;.
\end{align*}
$\\ $
Insbesondere stellen wir fest, dass für $0\leq k<P$ in der Ableitungsfolge
\[
u\left(x\right)^{\left(k\right)}=
\left(\left(x-\sqrt[M]{a}\right)^{{P}}\right)^{\left(k\right)}
=\frac{P\,!}{\left(P-k\right)\,!}\left(x-\sqrt[M]{a}\right)^{P-k}
\]
stets der Faktor
$\left(x-\sqrt[M]{a}\right)$ auftritt und somit für $x=\sqrt[M]{a}$ die Terme verschwinden. Nur im Falle $k=P$ haben wir einen von Null verschiedenen Wert, nämlich
\[
{\left.u{{\left(x\right)}^{\left(P\right)}}\right|}_{x=\sqrt[M]{a}}=\frac{P\,!}{\left(P-P\right)\,!}=P\,!\;.
\]
Von der Summe
\[
\left(f\left(x\right)^{P}\right)^{\left({P}\right)}
=\sum\limits_{k=0}^{{P}}\binom{{P}}{k}
\,\left(\left(x-\sqrt[M]{a}\right)^{P}\right)^{\left(k\right)}
\times
\left(\left(\frac{-M}{\sqrt[M]{a}}+\frac{f''{\left(\xi\right)}}{2}\left(x-\sqrt[M]{a}\right)\right)^{{P}}\right)^{\left({P-k}\right)}
\]
bleibt für $x=\sqrt[M]{a}$ somit nur der Term für $k=P$ übrig. Wir erhalten also
\begin{align*}
\left.\left(f\left(x\right)^{P}\right)^{\left({P}\right)}\right|_{x=\sqrt[M]{a}}
&=
\underbrace{\left.\left(\left(x-\sqrt[M]{a}\right)^{P}\right)^{\left(P\right)}\right|_{x=\sqrt[M]{a}}}_{=P\,!}
\times
\underbrace{\left.\left(\left(\frac{-M}{\sqrt[M]{a}}
+\frac{f''{\left(\xi\right)}}{2}\left(x-\sqrt[M]{a}\right)\right)^{{P}}\right)\right|_{x=\sqrt[M]{a}}}_{\left(\frac{-M}{\sqrt[M]{a}}\right)^P}
\\
&=P!\left(\frac{-M}{\sqrt[M]{a}}\right)^P=P\,!\frac{\left(-1\right)^P M^P}{\sqrt[M]{a^P}}\;.
\end{align*}
Unser ursprüngliches Anliegen war die Berechnung von
\[
F\left(\sqrt[M]{a}\right)^{\left(P+1\right)}
=\left.F\left(x\right)^{\left(P+1\right)}\right|_{x={\sqrt[M]{a}}}
=\left.F'\left(x\right)^{\left(P\right)}\right|_{x={\sqrt[M]{a}}}\;.
\]
Jetzt brauchen wir nur noch einzusetzen
\begin{align*}
F\left(\sqrt[M]{a}\right)^{\left(P+1\right)}
&=\left(\prod_{k=1}^{P}\left(1+{\frac{1}{k\,\cdot M}}\right)\cdot
\left.\left(1-{\frac{a}{{x}^{M}}}\right)^{P}\right)^{\left(P\right)}\right|_{x={\sqrt[M]{a}}}
\\
&=\prod_{k=1}^{P}\left(1+{\frac{1}{k\,\cdot M}}\right)\cdot P\,!\cdot\frac{\left(-1\right)^P\,M^P}{\sqrt[M]{a^P}}\;.
\end{align*}
Wir können einige Umformungen vornehmen. Es gilt nämlich
\[
\prod_{\ell=1}^{P}\left(\frac{1}{1+\ell\cdot M}\right)
=\prod_{\ell=1}^{P}\left(\frac{1+\ell\cdot M}{\ell\cdot M}\right)
=\frac{1}{P\,!\,M^P}\prod_{\ell=1}^{P}\left(1+\ell\cdot M\right)\;.
\]
Wir erhalten somit
\begin{align*}
F\left(\sqrt[M]{a}\right)^{\left(P+1\right)}
=\frac{1}{P!\,M^P}\prod_{\ell=1}^{P}\left(1+\ell\cdot M\right)\cdot\,P\,!\cdot\frac{\left(-1\right)^P\,M^P}{\sqrt[M]{a^P}}
=\frac{\left(-1\right)^P}{\sqrt[M]{a^P}}\,\prod_{\ell=1}^{P}\left(1+\ell\cdot M\right)\;.
\end{align*}
Fassen wir alle Ergebnisse zusammen, so ergibt sich das gewünschte Ergebnis
\begin{align*}
F{\left(\sqrt[M]{a}\right)}^{\left(k\right)}=\left\{
\begin{aligned}
0\qquad\qquad\qquad\qquad 1\,\le &\,k\,\leq{P}
\\
\\
\frac{\left(-1\right)^P}{\sqrt[M]{a^P}}\prod_{\ell=1}^{P}\left(1+\ell\cdot M\right)
\qquad\qquad\quad k=&P+1
\\ 
\end{aligned} 
\right.
\end{align*}
Damit ist der Beweis für die Ableitungen abgeschlossen.
\subsection{Konvergenz}
Den Beweis führen wir mit dem Fixpunktsatz von Banach \cite{HJ} und zwar in einer speziellen Form für stetig differenzierbare Abbildungen in $\mathbb{R}$.
$ \\ $
$ \\ $
\underline{Fixpunktsatz von Banach:}
\\
Es sei $U\subseteq R$ eine abgeschlossene Teilmenge. Weiter sei $F:U\rightarrow R$ eine Abbildung mit folgenden Eigenschaften
\[
F\left(U\right)\subseteq U \quad (Selbstabbildung)
\]
Es existiere ein $0<L<1$, so dass gilt
\[
\left|F\,'\left(x\right)\right|\leq L\; \text{für alle} \, x\in U \quad (Kontraktionseigenschaft)
\]
Dann gilt:
\begin{itemize}
\item[(a)]
Es gibt genau einen Fixpunkt ${{x}^{*}}\in U$ von $F$
\item[(b)]
Für jeden Startwert $x_0\in \mathbb{U}$ konvergiert die Folge
$x_{n+1}=F\left(x_{{n}}\right)$
gegen $x^*$.
\end{itemize}
Jetzt kommen wir zu unserem Anliegen zurück. Wir haben bereits berechnet $F\left(\sqrt[M]{a}\right)=\sqrt[M]{a}$ und $F\,'\left(\sqrt[M]{a}\right)=0$. Dann existiert aus Gründen der Stetigkeit um den Fixpunkt $\sqrt[M]{a}$ eine Umgebung
\[U_{\delta}=\left\{x\in\mathbb{R}:\left|x-\sqrt[M]{a}\right|\leq\delta\right\}\]
mit der Eigenschaft: $F\left(U\right)\subseteq U$ und $\left|F\,'\left(x\right)\right|<L<1$ für alle $x\in U_{\delta}$.
Somit sind die Voraussetzungen an den Fixpunktsatz erfüllt und hieraus folgt die Konvergenz.
\subsection{Konvergenzordnung}
Als letztes haben wir zu zeigen, dass die Konvergenzordnung genau $\left(P+1\right)$ beträgt.
Wir greifen auf bereits bewährtes zurück: Taylorreihen. Wir entwickeln $F\left(x\right)$ in eine Taylorreihe um den Fixpunkt  $x_0={\sqrt[M]{a}}$ wie folgt.
\begin{align*}
F\left(x\right)&=
F\left({\sqrt[M]{a}}\right)
+F'\left({\sqrt[M]{a}}\right)\left(x-{\sqrt[M]{a}}\right)
+\ldots
\\
&+F\left({\sqrt[M]{a}}\right)^{\left(P\right)}
\frac{\left(x-{\sqrt[M]{a}}\right)^{P}}{P\,!}
+F\left(\xi\right)^{\left(P+1\right)}
\frac{\left(x-{\sqrt[M]{a}}\right)^{P+1}}{\left(P+1\right)!}\;,
\end{align*}
mit einem $\xi$ zwischen $x$ und ${\sqrt[M]{a}}$.
Wir erinnern uns an die Fixpunkteigenschaft $F\left({\sqrt[M]{a}}\right)={\sqrt[M]{a}}$ und unsere Ableitungsergebnisse:
\begin{align*}
F{\left(\sqrt[M]{a}\right)}^{\left(k\right)}=\left\{
\begin{aligned}
0\qquad\qquad\qquad\qquad 0\,\le &\,k\,\leq{P}
\\
\\
\frac{\left(-1\right)^P}{\sqrt[M]{a^P}}\prod_{\ell=1}^{P}\left(1+\ell\cdot M\right)
\qquad\qquad\quad k=&P+1
\\ 
\end{aligned} 
\right.
\end{align*}
Die Taylorreihe
\begin{align*}
F\left(x\right)
&=
\underbrace{F\left({\sqrt[M]{a}}\right)}_{={\sqrt[M]{a}}}
+\underbrace{F'\left({\sqrt[M]{a}}\right)\left(x-{\sqrt[M]{a}}\right)
+\ldots
+F\left({\sqrt[M]{a}}\right)^{\left(P\right)}
\frac{\left(x-{\sqrt[M]{a}}\right)^{P}}{P\,!}}_{=0}
\\
&+F\left(\xi\right)^{\left(P+1\right)}
\frac{\left(x-{\sqrt[M]{a}}\right)^{P+1}}{\left(P+1\right)\,!}\,,
\end{align*}
vereinfacht sich damit zu
 \[
F\left(x\right)=
{\sqrt[M]{a}}
+F\left(\xi\right)^{\left(P+1\right)}
\frac{\left(x-{\sqrt[M]{a}}\right)^{P+1}}{\left(P+1\right)!}\;.
\]
Nun setzen wir $x=x_n$ hinreichend nahe bei ${\sqrt[M]{a}}$ ein und erhalten
\begin{align*}
\underbrace{F\left(x_n\right)}_{x_{n+1}}=
{\sqrt[M]{a}}
+F\left(\xi_n\right)^{\left(P+1\right)}
\frac{\left(x_n-{\sqrt[M]{a}}\right)^{P+1}}{\left(P+1\right)!}
\\
x_{n+1}={\sqrt[M]{a}}
+F\left(\xi_n\right)^{\left(P+1\right)}
\frac{\left(x_n-{\sqrt[M]{a}}\right)^{P+1}}{\left(P+1\right)!}\;.
\end{align*}
Wir bringen ${\sqrt[M]{a}}$ auf die linke Seite und dividieren durch
$\left(x_n-{\sqrt[M]{a}}\right)^{P+1}$. Das ergibt dann
\[
\frac{x_{n+1}-{\sqrt[M]{a}}}{\left(x_n-{\sqrt[M]{a}}\right)^{P+1}}=
F\left(\xi_n\right)^{\left(P+1\right)}
\frac{1}{\left(P+1\right)!}\;.
\]
Mit zunehmendem $n$ strebt $x_n\rightarrow {\sqrt[M]{a}}$ und da $\xi_n$ zwischen $x_n$ und ${\sqrt[M]{a}}$ eingesperrt ist, muss $\xi_n$ notgedrungen ebenfalls gegen ${\sqrt[M]{a}}$ streben. Wir erhalten also
\begin{align*}
\underset{n\to \infty}{\mathop{\lim }}\; \frac{x_{n+1}-{\sqrt[M]{a}}}{\left(x_n-{\sqrt[M]{a}}\right)^{P+1}}
=\underset{n\to \infty}{\mathop{\lim }}\;
F\left(\xi_n\right)^{\left(P+1\right)}
\frac{1}{\left(P+1\right)!}
=F\left({\sqrt[M]{a}}\right)^{\left(P+1\right)}
\frac{1}{\left(P+1\right)!}\;.
\end{align*}
Nun haben wir bereits berechnet
\[
F{\left({\sqrt[M]{a}}\right)}^{\left({P+1}\right)}
=\frac{\left(-1\right)^P}{\sqrt[M]{a^P}}\prod_{\ell=1}^{P}\left(1+\ell\cdot M\right)\;.
\]
Dies eingesetzt ergibt dann
\begin{align*}
\underset{n\to \infty}{\mathop{\lim }}\;\;
\frac{x_{n+1}-{\sqrt[M]{a}}}{\left(x_n-{\sqrt[M]{a}}\right)^{P+1}}
=\frac{\left(-1\right)^P}{\left(P+1\right)!}\frac{1}{\sqrt[M]{a^P}}\prod_{\ell=1}^{P}\left(1+\ell\cdot M\right)\;.
\end{align*}
Hieraus folgt unmittelbar die Konvergenzordnung $\left(P+1\right)$. Da $F{\left({\sqrt[M]{a}}\right)}^{\left({P+1}\right)}\neq 0$ haben wir genau Konvergenzordnung $\left(P+1\right)$. Der auf der rechten Seite stehende konstante Term
\[
\frac{1}{\left(P+1\right)\,!}\,\frac{\left(-1\right)^P}{\sqrt[M]{a^P}}\prod_{\ell=1}^{P}\left(1+\ell\cdot M\right)
\]
wird als asymtotische Konvergenzrate bezeichnet.
\section{Praktische Anwendung des Verfahrens}
Es seien $a>0$ und $M \in \mathbb{N}$ vorgegeben. Wir wollen $\sqrt[M]{a}$ berechnen. Zusätzlich geben wir $P\in\mathbb{N}$ vor. Da ergibt dann die Konvergenzordnung $\left(P+1\right)$.
\subsection{Übersicht über die Polynome bis Konvergenzordnung 5}
Bevor wir mit der Praxis beginnen verschaffen wir uns zunächst einen Überblick über die Polynome, die sich für $P \in \{1,2,3,4\}$ ergeben. Wir erhalten:
\begin{center}
\text{\underline{quadratisch}}
\end{center}
\[
\left(1+\frac{1}{M}\right)\left( x-\frac{x^{M+1}}{a\cdot M}\right)\,;
\]

\begin{center}
\text{\underline{kubisch}}
\end{center}
\[
\left(1+\frac{1}{M}\right)\left(1+\frac{1}{2M}\right)\left(x-\frac{2\,x^{M+1}}{\left(M+1\right)a}+\frac{x^{2M+1}}{\left(2\,M+1\right)a^{2}}\right)\,;
\]

\begin{center}
\text{\underline{biquadratisch}}
\end{center}
\[
\left(1+\frac{1}{M}\right)\left(1+\frac{1}{2\,M}\right)\left(1+\frac{1}{3\,M}\right)\left(x-\frac{3\,x^{M+1}}{\left(M+1\right)a}+\frac{3\,x^{2\,M+1}}{\left(2\,M+1\right)a^{2}}-\frac{x^{3\,M+1}}{\left(3\,M+1\right)a^{3}}\right)\,;
\]

\begin{center}
\text{\underline{quintisch}}
\end{center}
\[
\left(1+\frac{1}{M}\right)\left(1+\frac{1}{2\,M}\right)\left(1+\frac{1}{3\,M}\right)\left(1+\frac{1}{4\,M}\right)\times
\]
\[
\left(x-\frac{4x^{M+1}}{\left(M+1\right)a}+\frac{6\,x^{2\,M+1}}{\left(2\,M+1\right)a^{2}}-\frac{4x^{3\,M+1}}{\left(3\,M+1\right)a^{3}}+\frac{x^{4\,M+1}}{\left(4\,M+1\right)a^{4}}\right)\,;
\]
\subsection{Praktische Vorgehensweise allgemein}
Es liegen vor: $a>0$ sowie $M\in{\mathbb{N}}$ und $P\in{\mathbb{N}}$.
$ \\ $
Mit diesen Werten definieren wir zuerst die Funktion
\[
f\left(t\right)=1-\frac{t^M}{a}\;.
\]
und konstruieren die Fixpunkt-Funktion
\[
F\left(x\right)
=\prod_{\ell=1}^{P}\left(1+\frac{1}{\ell\cdot M}\right)
\int\limits_{0}^{x}\!f\left(t\right)^{P}{\rm d}t
=\prod_{\ell=1}^{P}\left(1+\frac{1}{\ell\cdot M}\right)
\int\limits_{0}^{x}\!\left(1-{\frac{{t}^{M}}{a}}\right)^{P}{\rm d}t\;.
\]
Dazu wenden wir auf
\[
\left(1-\frac{t^M}{a}\right)^{P}
\]
den binomischen Lehrsatz an und erhalten die Summe
\begin{align*}
\left(1-\frac{t^M}{a}\right)^{P}
&=\sum\limits_{k=0}^{P}\left(-1\right)^k\,\cdot\binom{P}{k}\,\left(\frac{t^M}{a}\right)^k
\\
&=\sum\limits_{k=0}^{P}\left(-1\right)^k\,\cdot\binom{P}{k}\cdot\frac{t^{k\,\cdot M}}{a^k}
\\
&=\sum\limits_{k=0}^{P}\frac{\left(-1\right)^k}{a^k}\cdot\binom{P}{k}\cdot t^{k\,\cdot M}\;.
\end{align*}
Die Summe lässt sich leicht integrieren
\[
\int\limits_{0}^{x}\!\left(1-{\frac{{t}^{M}}{a}}\right)^{P}{\rm d}t
=\int\limits_{0}^{x}\!\sum\limits_{k=0}^{P}\frac{\left(-1\right)^k}{a^k}\cdot\binom{P}{k}\cdot t^{\,k\,\cdot M}{\rm d}t
=\sum\limits_{k=0}^{P}\frac{\left(-1\right)^k}{a^k}\cdot\binom{P}{k}\cdot\frac{x^{\,k\,\cdot M+1}}{k\,\cdot M+1}\;.
\]
Jetzt fehlt noch die Multiplikation mit dem Produkt
\[
\prod_{\ell=1}^{P}\left(1+\frac{1}{\ell\cdot M}\right)
\]
und die Fixpunktfunktion ist fertig:
\[
F\left(x\right)=\prod_{\ell=1}^{P}\left(1+\frac{1}{\ell\cdot M}\right)
\sum\limits_{k=0}^{P}\frac{\left(-1\right)^k}{a^k}\cdot\binom{P}{k}\cdot\frac{x^{\,k\,\cdot M+1}}{k\,\cdot M+1}\;.
\]
Wir stellen wir fest, dass sich die etwas kompliziert gebauten Koeffizienten während der gesamten Berechnung nicht ändern. Es ist deshalb zweckmäßig diese vorab zu berechnen und die einzelnen Elemente in einem Koeffizientenvektor abzuspeichern:
\[
c_{\,k}=\prod_{\ell=1}^{P}\left(1+\frac{1}{\ell\cdot M}\right)
\frac{\left(-1\right)^k}{a^k}\cdot\binom{P}{k}\cdot\frac{1}{k\,\cdot M+1}\;;\;k=0,1,\ldots,P\]
Die Fixpunktfunktion vereinfacht sich damit zu
\[F\left(x\right)=\sum\limits_{k=0}^{P}c_{\,k}\cdot x^{\,k\,\cdot M+1}\;.\]
Das ist ein Polynom vom Grad $\left(P\cdot M+1\right)$ und besteht aus $\left(P+1\right)$ Summanden.
$ \\ $
$ \\ $
\underline{Hinweis:} Die Berechnung von $F\left(x\right)$ ist somit eine Polynomauswertung. Für Polynomauswertungen gibt es effiziente Methoden (Horner-Schema).
$ \\ $
$ \\ $
Bei der praktischen Berechnung ist eine Indexverschiebung sinnvoll. Anstatt
\[x_{n+1}=\sum\limits_{k=0}^{P}c_{\,k}\cdot x_{n}^{\,k\,\cdot M+1}\;;\;n=0,1,\ldots\]
verwenden wir die Schreibweise
\[x_{n}=\sum\limits_{k=0}^{P}c_{\,k}\cdot x_{n-1}^{\,k\,\cdot M+1}\;;\;n=1,2,\ldots\]
Damit läuft der Index von $x_n$ synchron mit der Anzahl der Iterationen.
$ \\ $
$ \\ $
Die Ausführung ist denkbar einfach.
$ \\ $
Wir geben einen geeigneten Startwert $x_0$ vor und berechnen

\[\underline{step\;\;1}\]
\[x_{1}=\sum\limits_{k=0}^{P}c_{\,k}\cdot x_{0}^{\;k\,\cdot M+1}\;;\]

\[\underline{step\;\;2}\]
\[x_{2}=\sum\limits_{k=0}^{P}c_{\,k}\cdot x_{1}^{\;k\,\cdot M+1}\;;\]

\[\underline{step\;\;3}\]
\[x_{3}=\sum\limits_{k=0}^{P}c_{\,k}\cdot x_{2}^{\;k\,\cdot M+1}\;;\]


\[\vdots\]

\[\underline{step\;\;n}\]
\[x_{n}=\sum\limits_{k=0}^{P}c_{\,k}\cdot x_{n-1}^{\;k\,\cdot M+1}\;;\]
$ \\ $
Mit zunehmender Anzahl $n$ der Iterationsschritte strebt $x_{n+1}\rightarrow\sqrt[M]a$ mit Konvergenzordnung $\left(P+1\right)$. Bei beliebig langer Rechenzeit und beliebig viel Speicherplatz können wir damit $\sqrt[M]a$ beliebig genau berechnen. Leider steht uns nichts dergleichen zur Verfügung. Wir müssen uns mit einer endlichen Anzahl von Iterationen und einem endlichen Wert für $\sqrt[M]a$ begnügen. Dazu geben wir $\epsilon>0$ vor und führen die Iteration solange aus bis erstmals gilt:
\[\left|x_{n}-x_{n-1}\right|<\epsilon\]
\subsection{einzelne Wurzelbrechnungen}
\textbf{Beispiel $1\,:$ Berechnung von $\sqrt[3]{10}$}
$ \\ $
Im ersten Beispiel berechnen wir $\sqrt[3]{10}$. Es ist also $a=10$ und $M=3$. We geben $P=1$ vor. Das führt zur Konvergenzordnung $\left(P+1\right)=2$ (quadratisch). Jeder zusätzlicher Iterationschritt verdoppelt näherungsweise die Anzahl gültiger Stellen. Einsetzen von $a=10\;,M=3$ ergibt zunächst $f\left(t\right)=1-\frac{t^3}{10}\;.$ Mit $P=1$ ist das Fixpunktpolynom bestimmt als:
\[F\left(x\right)
=\prod_{k=1}^{1}\left(1+\frac{1}{3\cdot k}\right)
\int\limits_{0}^{x}\!\left(1-{\frac{{t}^{3}}{10}}\right)
=\left(1+\frac{1}{3}\right)\left(x-\frac{1}{40}\,x^4\right)
=\frac{4}{3}\,x-\frac{1}{30}\,x^4\;.\]
Wir sehen uns die beteiligten Funktionen zusammen mit der Winkelhalbierenden $y=x$ an:
\begin{center}
\includegraphics[width=0.95\linewidth]{"fp_integral_sqrt10_M3_P1.eps"}
\end{center}
\begin{center}
\textbf{\underline{Abbildung 1}}
\end{center}
Wir wählen als Startwert $x_0=2$ und berechnen die ersten sechs Iterationen:

\[\underline{step\;\;1}\]
\[x_{1}=\frac{4}{3}x_{0}-\frac{1}{30}x_{0}^{4}\]
\[x_{1}=2.133333333333333333333333333333333333333\]
\[{|x_{1}-x_{0}|}=0.1333333333333333333333333333333333333333\]
\[\underline{step\;\;2}\]
\[x_{2}=\frac{4}{3}x_{1}-\frac{1}{30}x_{1}^{4}\]
\[x_{2}=2.154024032921810699588477366255144032922\]
\[{|x_{2}-x_{1}|}=0.02069069958847736625514403292181069958848\]
\[\underline{step\;\;3}\]
\[x_{3}=\frac{4}{3}x_{2}-\frac{1}{30}x_{2}^{4}\]
\[x_{3}=2.154434533500953092649669501763572523986\]
\[{|x_{3}-x_{2}|}=0.0004105005791423930611921355084284910642133\]
\[\underline{step\;\;4}\]
\[x_{4}=\frac{4}{3}x_{3}-\frac{1}{30}x_{3}^{4}\]
\[x_{4}=2.154434690031860976181374509716973801410\]
\[{|x_{4}-x_{3}|}={1.565309078835317050079534012774237318926\times10^{-7}}\]
\[\underline{step\;\;5}\]
\[x_{5}=\frac{4}{3}x_{4}-\frac{1}{30}x_{4}^{4}\]
\[x_{5}=2.154434690031883721759293566039074794849\]
\[{|x_{5}-x_{4}|}={2.274557791905632210099343907978738060749\times10^{-14}}\]
\[\underline{step\;\;6}\]
\[x_{6}=\frac{4}{3}x_{5}-\frac{1}{30}x_{5}^{4}\]
\[x_{6}=2.154434690031883721759293566519350495259\]
\[{|x_{6}-x_{5}|}={4.802757004105093077094334087308664908888\times10^{-28}}\]
$ \\ $
An den Differenzen $\left|x_5-x_4\right|$ und $\left|x_6-x_5\right|$ können wir sehr schön die Konvergenzordnung $2$ ablesen. Der negative Exponent verdoppelt sich mit jedem Iterationsschritt. 
$ \\ $
$ \\ $
$ \\ $
\textbf{Beispiel $2\,:$ Berechnung einer Million Stellen für $\sqrt{2}$}
$ \\ $
Im nächsten Beispiel wollen wir $\sqrt{2}$ auf mindestens eine Million Dezimalstellen genau berechnen. Das bedeutet: $\epsilon=10^{-1000000}$. Es ist $a=2$ und $M=2$. Das ergibt $f\left(t\right)=1-\frac{t^2}{2}.$ Wir geben  $P=3$ vor und erhalten Konvergenzordnung $3+1=4$. Beim Startwert gehen wir diesmal viel näher an $\sqrt{2}$ ran. Wir wählen den auf 16 Dezimalstellen genauen Näherungswert $x_0=1.414213562373095$. Da sich die Anzahl der gültigen Stellen mit jedem Iterationsschritt im Prinzip vervierfacht, hoffen wir nach 9 Iterationen am Ziel zu sein. Wir beginnen mit der Konstuktion der Fixpunktfunktion
\[
F\left(x\right)
=\prod_{\ell=1}^{3}\left(1+{\frac{1}{2\,\ell}}\right)\int\limits_{0}^{x}\!\left(1-{\frac{{t}^{2}}{2}}\right)^{3}\,{\rm d}t
=\ldots=\frac{35}{16}\,x-\frac{35}{32}\,{x}^{3}+\frac{21}{64}\,{x}^{5}-\frac{5}{128}\,{x}^{7}\;.
\]
Wir sehen uns die beteiligten Funktionen in einer weiteren Abbildung an.
\begin{center}
\includegraphics[width=0.95\linewidth]{"fp_integral_sqrt02_M2_P3.eps"}
\textbf{\underline{Abbildung 2}}
\end{center}
\[\]
\[\]
Mit Startwert $x_0=1.414213562373095$ erhalten wir die Iterationsfolge:
\[\underline{step\;\;1}\]
\[x_{1}=\frac{35}{16}x_{0}-\frac{35}{32}x_{0}^{3}+\frac{21}{64}x_{0}^{5}-\frac{5}{128}x_{0}^{7}\]
\[x_{1}=1.414213562373095048801688724209698078570\]
\[{|x_{1}-x_{0}|}={4.880168872420969807856967187537694807318\times10^{-17}}\]
\[\underline{step\;\;2}\]
\[x_{2}=\frac{35}{16}x_{1}-\frac{35}{32}x_{1}^{3}+\frac{21}{64}x_{1}^{5}-\frac{5}{128}x_{1}^{7}\]
\[x_{2}=1.414213562373095048801688724209698078570\]
\[{|x_{2}-x_{1}|}={8.773491625654111352087407579690431191435\times10^{-66}}\]
\[\underline{step\;\;3}\]
\[x_{3}=\frac{35}{16}x_{2}-\frac{35}{32}x_{2}^{3}+\frac{21}{64}x_{2}^{5}-\frac{5}{128}x_{2}^{7}\]
\[x_{3}=1.414213562373095048801688724209698078570\]
\[{|x_{3}-x_{2}|}={9.164798637556653681657805406878049888878\times10^{-261}}\]
\[\underline{step\;\;4}\]
\[x_{4}=\frac{35}{16}x_{3}-\frac{35}{32}x_{3}^{3}+\frac{21}{64}x_{3}^{5}-\frac{5}{128}x_{3}^{7}\]
\[x_{4}=1.414213562373095048801688724209698078570\]
\[{|x_{4}-x_{3}|}={1.091251298365935101705686744387078883102\times10^{-1040}}\]

\[\underline{step\;\;5}\]
\[x_{5}=\frac{35}{16}x_{4}-\frac{35}{32}x_{4}^{3}+\frac{21}{64}x_{4}^{5}-\frac{5}{128}x_{4}^{7}\]
\[x_{5}=1.414213562373095048801688724209698078570\]
\[{|x_{5}-x_{4}|}={2.193472316487722705810599621121648551289\times10^{-4160}}\]
\[\underline{step\;\;6}\]
\[x_{6}=\frac{35}{16}x_{5}-\frac{35}{32}x_{5}^{3}+\frac{21}{64}x_{5}^{5}-\frac{5}{128}x_{5}^{7}\]
\[x_{6}=1.414213562373095048801688724209698078570\]
\[{|x_{6}-x_{5}|}={3.580648536099876136173035995717511426715\times10^{-16639}}\]

\[$$ $$ \underline{step\;\;7}\]
\[x_{7}=\frac{35}{16}x_{6}-\frac{35}{32}x_{6}^{3}+\frac{21}{64}x_{6}^{5}-\frac{5}{128}x_{6}^{7}\]
\[x_{7}=1.414213562373095048801688724209698078570\]
\[{|x_{7}-x_{6}|}={2.542610528450840832485991523758935060375\times10^{-66554}}\]
\[\underline{step\;\;8}\]
\[x_{8}=\frac{35}{16}x_{7}-\frac{35}{32}x_{7}^{3}+\frac{21}{64}x_{7}^{5}-\frac{5}{128}x_{7}^{7}\]
\[x_{8}=1.414213562373095048801688724209698078570\]
\[{|x_{8}-x_{7}|}={6.464760315447686077979797373449536529093\times10^{-266215}}\]
\[\underline{step\;\;9}\]
\[x_{9}=\frac{35}{16}x_{8}-\frac{35}{32}x_{8}^{3}+\frac{21}{64}x_{8}^{5}-\frac{5}{128}x_{8}^{7}\]
\[x_{9}=1.414213562373095048801688724209698078570\]
\[{|x_{9}-x_{8}|}={2.701735162639912537134047073288055961734\times10^{-1064857}}\]
$\\ $
Das Ergebnis liegt nach 3 Sekunden vor. Mit 9 Iterationen haben wir tatsächlich über eine Million Stellen für $\sqrt{2}$ erhalten. An den negativen Exponenten der Differenzen können wir die Konvergenzordnung 4 ablesen. Die negativen Exponenten vervierfachen sich prinzipiell mit jedem Schritt.
$ \\ $
$ \\ $
Alle Berechnungen wurden mit einem Eigenbau-PC mit folgender Hardware ausgeführt: Motherboard ASUS PRIME A320M-K mit CPU AMD Ryzen 5 5600G 6 CORE 3.90-4.40 GHz und 32 GB RAM. Die verwendete Software war MAPLE 2025.2 von  Maplesoft, Waterloo Maple Inc. Die Software wurde ebenfalls für die Konvertierung der mathematischen Terme nach LaTeX verwendet. 


\begin{thebibliography}{99}

\bibitem{WG} Walter Gautschi: Numerical Analysis Second Edition, {\it Birkhäuser 2012}

\bibitem{AS} Milton Abramowitz and Irene Stegun: Handbook of Mathematical Functions with Formulas, Graphs, and Mathematical Tables, {\it National Institute of Standards and Technology 1965}

\bibitem{EA} Emil Artin: The Gamma Function, {\it Holt, Rinehart and Winston, Inc. 1965}

\bibitem{HJ} HUGO D. JUNGHENN: A COURSE IN REAL ANALYSIS, {\it CRC Press Taylor and Francis Group 2015}

\bibitem{formula} Schaums's Outline Series: Mathematical Handbook of Formulas and Tables, {\it McGraw-Hill Book Company 1968}

\bibitem{JM} J. M. McDonough: LECTURES IN BASIC COMPUTATIONAL NUMERICAL ANALYSIS, {\it Departments of Mechanical Engineering and Mathematics University of Kentucky 2007}

\end{thebibliography}
\end{document}